\documentclass{bmcart}

\usepackage{amsthm, amsmath, amssymb}
\usepackage[per-mode=symbol]{siunitx}
\usepackage{booktabs}
\usepackage[noline, linesnumbered, noend, ruled]{algorithm2e}
\usepackage{subcaption}
\usepackage{graphicx}
\usepackage{pgfplots}
\usepackage{nicefrac}
\usepackage{cite}
\usepackage{braket}

\RequirePackage{hyperref}
\usepackage{hyperref}
\usepackage[utf8]{inputenc} 

\startlocaldefs

\newcommand{\dx}[1]{\left.\text{d}#1 \right.}

\newcommand{\qe}[1]{``#1''}
\newcommand{\R}{\mathbb{R}}
\newcommand{\normal}{n}
\newcommand{\abs}[1]{\left\lvert #1 \right\rvert}
\newcommand{\norm}[2]{\left\lvert\left\lvert #1 \right\rvert\right\rvert_{#2}}

\newcommand{\subcool}{\text{cool}}
\newcommand{\subrad}{\text{rad}}
\newcommand{\subamb}{\text{amb}}
\newcommand{\subblood}{\text{b}}
\newcommand{\subnative}{n}
\newcommand{\subcoag}{c}
\newcommand{\subext}{\text{ext}}
\newcommand{\subdes}{\text{meas}}
\newcommand{\subend}{\text{end}}

\newcommand{\temperature}{T}
\newcommand{\radiation}{\varphi}
\newcommand{\damage}{\omega}

\newcommand{\adtemp}{p}
\newcommand{\adrad}{\psi}
\newcommand{\dmgfun}{\delta}

\newcommand{\conductivity}{\kappa}
\newcommand{\density}{\rho}
\newcommand{\heatcapacity}{C_p}
\newcommand{\perfusion}{\xi}

\newcommand{\absorption}[1]{\mu_{a#1}}
\newcommand{\scattering}[1]{\mu_{s#1}}
\newcommand{\anisotropy}[1]{g_{#1}}
\newcommand{\timehorizon}{\tau}
\newcommand{\htc}{\alpha}
\newcommand{\diffusivity}{D}
\newcommand{\laserpower}{q_{\text{app}}}
\newcommand{\coolingfactor}{\beta_q}
\newcommand{\laseron}{t_\text{on}}
\newcommand{\laseroff}{t_\text{off}}
\newcommand{\frequencyfactor}{A}
\newcommand{\activationenergy}{E_a}
\newcommand{\gasconstant}{R}

\newenvironment{eqsystem}[1][]
{\begin{equation}
 \label{#1}
 \left\lbrace
 \begin{alignedat}{2}
}
{\end{alignedat}
\right.
\end{equation}
}

\newenvironment{eqsystem*}
{\begin{equation*}
	\left\lbrace
	\begin{alignedat}{2}
}
{\end{alignedat}
	\right.
	\end{equation*}
}

\newtheorem{Theorem}{Theorem}[section]
\newtheorem{Proposition}[Theorem]{Proposition}

\endlocaldefs

\begin{document}

\begin{frontmatter}

\begin{fmbox}
\dochead{Research}


\title{Identification of the Blood Perfusion Rate for Laser-Induced Thermotherapy in the Liver}

\author[
	addressref={tuk}, 
	email={andres@mathematik.uni-kl.de}
]{\inits{M}\fnm{Matthias} \snm{Andres}}
\author[
	addressref={itwm, tuk},
	email={sebastian.blauth@itwm.fraunhofer.de}
]{\inits{S}\fnm{Sebastian} \snm{Blauth}}
\author[
	addressref={itwm},
	email={christian.leithäuser@itwm.fraunhofer.de}
]{\inits{C}\fnm{Christian} \snm{Leithäuser}}
\author[
	addressref={itwm},
	corref={itwm},
	email={norbert.siedow@itwm.fraunhofer.de}
]{\inits{N}\fnm{Norbert} \snm{Siedow}}


\address[id=tuk]{%
  \orgname{TU Kaiserslautern},
  \street{Gottlieb-Daimler Straße 48},
  \postcode{67663}
  \city{Kaiserslautern},
  \cny{Germany}
}
\address[id=itwm]{
	\orgname{Fraunhofer Institute for Industrial Mathematics ITWM}, 
	\street{Fraunhofer Platz 1},                     %
	\postcode{67663}                                
	\city{Kaiserslautern},                              
	\cny{Germany}                                    
}


\begin{artnotes}
\end{artnotes}

\end{fmbox}


\begin{abstractbox}

\begin{abstract} 
Using PDE-constrained optimization we introduce a parameter identification approach which can identify the blood perfusion rate from MR thermometry data obtained during the treatment with laser-induced thermotherapy (LITT). The blood perfusion rate, i.e., the cooling effect induced by blood vessels, can be identified during the first stage of the treatment. This information can then be used by a simulation to monitor and predict the ongoing treatment. The approach is tested with synthetic measurements with and without artificial noise as input data.

\end{abstract}




\end{abstractbox}
%

\end{frontmatter}




\section{Introduction}
\label{sec:introduction}

Laser-induced interstitial thermotherapy (LITT) is a minimally invasive, local therapy used to destroy tumors through thermal ablation. For this, laser radiation is transmitted by an optical fiber to an application system that is inserted into the tumorous tissue. Absorption of the radiation by the tissue results in a temperature increase around the applicator which destroys the tumor cells due to coagulative effects. The goal of the therapy is to completely destroy the tumor while protecting the healthy tissue. To reach this goal, computer simulations can assist physicians in the planing and monitoring of treatments. However, such simulations can only yield reliable results if all necessary parameters are known. While typically good measurements are available for many of the tissue parameters, a critical role is the determination of the blood perfusion rate that models the cooling effect induced by blood vessels.

The blood perfusion rate is a nonlocal whose magnitude depends on the presence of blood vessels. Further, it depends on the shape and size of the vessels. The induced cooling effect significantly influences the temperature. The knowledge of the location of the blood vessels in the vicinity of tumorous tissue (and, thus, close to the applicator) is crucial for the performance of the therapy as well as for the reliability of a simulation model as e.g. discussed in \cite{perfusion_parameter, mohammed_verhey, shao2017computational, kroger2006numerical, shibib2017effect}. Unfortunately, the location relative to the applicator varies for each patient and treatment. The location of major vessels can be determined a-priori, e.g., with the help of image decomposition techniques. However, this is tedious and may give a rather bad approximation of the actual perfusion rate. 

A more promising approach for the identification of the blood perfusion rate is to use temperature measurements obtained by magnetic resonance (MR) thermometry. MR thermometry methods are based on MR measured parameters depending on temperature like the longitudinal relaxation time, the diffusion coefficient, or the proton resonance frequency (PRF) of tissue water. The linear temperature dependence of the proton resonance frequency and its near-independence with respect to tissue type make the PRF-based methods the preferred choice for many applications. For a more deeper understanding to MR-Thermometry we refer to the review paper \cite{Senne}.


The idea proposed in this paper is to identify the perfusion rate in a short time period during the beginning of the treatment from MR thermometry data. This information can then be used to simulate the remaining treatment. This is of great benefit for LITT as it can be integrated into an online therapy monitoring and prediction tool, individualized for every patient. In the following we derive a parameter identification approach for the blood perfusion rate from MR thermometry data using techniques from optimal control for partial differential equations. Synthetic measurements are used to demonstrate the approach.

This paper is structured as follows. In Section~\ref{sec:mathematical_model} we introduce our mathematical model of LITT. Section~\ref{sec:parameter_identification} gives the details of the parameter identification problem and its formulation as a PDE constrained optimization problem. The numerical solution of the problem is described in Section~\ref{sec:numerics}. The validation of our method is done in Section~\ref{sec:results}, where we discuss the results of some model problems using synthetic measurements.

\begin{figure}[!t]
	\centering
	\includegraphics[width=0.95\textwidth]{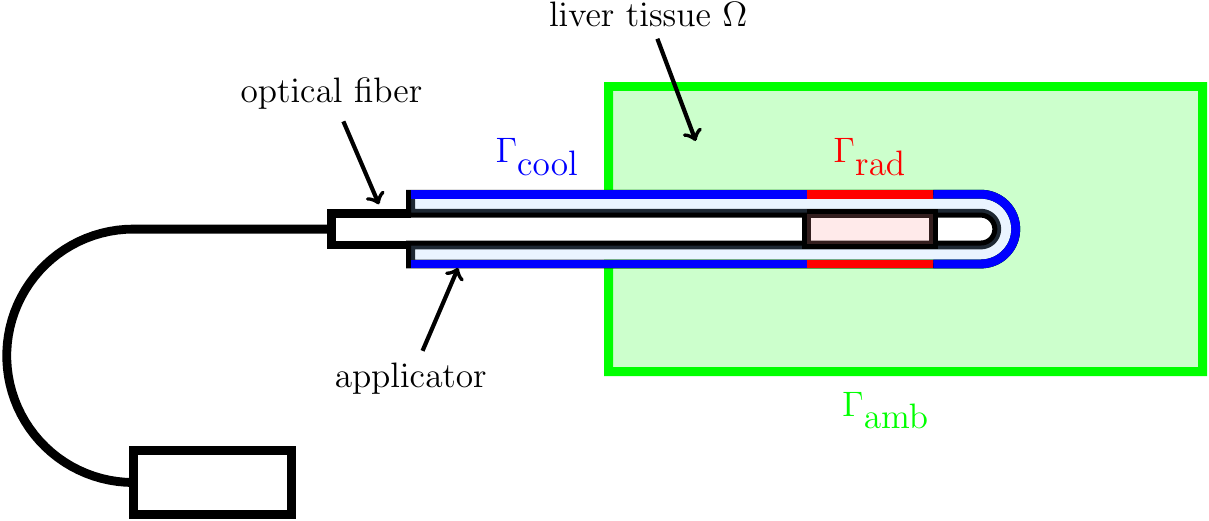}
	\caption{Schematic of the setup and the boundary decomposition.}
\end{figure}

\section{Mathematical Model}
\label{sec:mathematical_model}

For the modeling of the LITT we use the same model as \cite{huebner} which was proposed in \cite{fasano}. The model is summarized in the following. The liver tissue is denoted by $\Omega \subset \R^3$. Its boundary $\Gamma$ consists of the following three parts: The ambient boundary $\Gamma_\subamb$, i.e., the surface of the liver, as well as two parts corresponding to the interface between tissue and laser applicator (the latter is not part of $\Omega$): $\Gamma_\subcool$, the part of the applicator boundary that is cooled, and $\Gamma_\subrad$, the part of the boundary where radiation is emitted. Further, $t$ and $\tau$ denote the time [\si{\second}] and the time horizon of the simulation [\si{\second}], respectively. 

To model the tissue temperature during LITT we use Pennes' bio-heat equation \cite{pennes} that is coupled with a $P1$-approximation for the laser radiation \cite{niemz} and the Arrhenius law that models tissue damage \cite{fasano}. The bio-heat equation \cite{pennes} reads
\begin{eqsystem}[eq:bioheat]
	\density \heatcapacity \dot{\temperature} - \nabla \cdot (\conductivity \nabla \temperature) + \perfusion (\temperature - \temperature_\subblood) &= \absorption{} \radiation \quad &&\text{ in } (0, \timehorizon) \times \Omega, \\
	\conductivity \partial_\normal T + \htc (\temperature - \temperature_\subext) &= 0 \quad &&\text{ on } (0, \timehorizon) \times \Gamma, \\
	\temperature(0, \cdot) &= \temperature_0 \quad &&\text{ in } \Omega,
\end{eqsystem}
where $\temperature$ is the tissue temperature [\si{\kelvin}], and $\conductivity$, $\density$ and $\heatcapacity$ denote the liver's thermal conductivity [\si{\watt \per \meter \per \kelvin}], density [\si{\kilogram \per \cubic\meter}] and specific heat capacity [\si{\joule \per \kilogram \per \kelvin}], respectively. The term $\perfusion(\temperature - \temperature_\subblood)$ models the heat transfer between blood vessels and liver tissue, where $\perfusion$ denotes the perfusion rate [\si{\watt \per \kelvin \per \cubic\meter}] and $\temperature_\subblood$ is the blood temperature [\si{\kelvin}]. As the perfusion rate is typically unknown, we propose a method for identifying this quantity later on. The energy generated by the laser enters the equation as a source term, where $\absorption{}$ is the absorption coefficient [\si{1 \per \meter}] and $\radiation$ is the radiative energy [\si{\watt \per \square\meter}], which we explain in the following paragraph (cf.~\eqref{eq:radiation}). The term $\partial_\normal \temperature$ denotes the normal derivative of $\temperature$ in the direction of the outer unit normal vector $\normal$, i.e., $\partial_\normal \temperature = \normal \cdot \nabla \temperature$. The heat transfer coefficient [\si{\watt \per \kelvin \per \square\meter}] is denoted by $\htc$ and the external temperature [\si{\kelvin}] by $\temperature_\subext$. Note that these two parameters vary over the boundary, in particular, we have that $\htc = \htc_\subcool$ and $\temperature_\subext = \temperature_\subcool$ on $\Gamma_\subcool \cup \Gamma_\subrad$ as well as $\htc = \htc_\subamb$ and $\temperature_\subext = \temperature_\subamb$ on $\Gamma_\subamb$. Finally, the initial temperature of the tissue [\si{\kelvin}] is denoted by $\temperature_0$.

We model the radiative energy $\radiation$ via the $P1$-approximation \cite{niemz, modest}:
\begin{eqsystem}[eq:radiation]
	- \nabla \cdot (\diffusivity \nabla \radiation) + \absorption{} \radiation &= 0 \quad &&\text{ in } (0, \timehorizon) \times \Omega, \\
	D \partial_\normal \radiation &= q \quad &&\text{ on } (0, \timehorizon) \times \left( \Gamma_\subcool \cup \Gamma_\subrad \right), \\
	D \partial_\normal \radiation + \frac{1}{2} \radiation &= 0 \quad &&\text{ on } (0, \timehorizon) \times \Gamma_\subamb,
\end{eqsystem}
where $\diffusivity$ denotes the diffusion coefficient [\si{m}] that is defined as
\begin{equation*}
	D = \frac{1}{3\left( \absorption{} + \scattering{} \left( 1 - \anisotropy{} \right) \right)}.
\end{equation*}
Here, $\scattering{}$ is the scattering coefficient [\si{1 \per \meter}] and $\anisotropy{}$ is the anisotropy factor.
Furthermore, $q$ is a boundary source that models the radiation coming from the laser applicator. On $\Gamma_\subcool$ there is no radiation, so $q = 0$, and on $\Gamma_\subrad$ we have $q = \frac{\laserpower}{\abs{\Gamma_\subrad}}$, where $\abs{\Gamma_\subrad}$ denotes the area of the radiating surface $\Gamma_\subrad$ and $\laserpower$ is the effective laser power. The latter is given by
\begin{equation*}
	\laserpower(t) = 
	\begin{cases}
		(1 - \coolingfactor) \hat{q}_\text{app} \quad &\text{ if } \laseron \leq t \leq \laseroff,\\
		0 \quad &\text{ otherwise},
	\end{cases}
\end{equation*}
where $\hat{q}_\text{app}$ is the total power emitted by the laser [\si{\watt}] and $\coolingfactor$ is the coolant absorption factor that models the absorption of energy by the coolant (cf. \cite{huebner}). Additionally, $\laseron$ and $\laseroff$ denote the times where the laser is switched on and off, respectively.

\begin{table}[!b]
	\begin{tabular}{l r}
		\toprule
		parameter [unit]  & value\\
		\midrule
		general & \\
		\midrule
		tissue density $\density$ [\si{\kilogram \per \cubic\meter}] & \num{1.08e3} \\
		universal gas constant $\gasconstant$ [\si{\joule \per \mole \per \kelvin}] & \num{8.31} \\
		total laser power $\hat{q}_\text{app}$ [\si{\watt}] & \num{22} \\
		$\laseron$ [\si{\second}] & \num{25} \\
		$\laseroff$ [\si{\second}] & 1175 \\
		end of treatment $\timehorizon_\subend$ [\si{\second}] & 1200 \\
		perfusion rate in a blood vessel $\perfusion_\text{max}$ [\si{\watt \per \kelvin \per \cubic\meter}] & \num{6e4} \\
		\midrule
		thermal & \\
		\midrule
		heat conductivity $\conductivity$ [\si{\watt \per \meter \per \kelvin}] & \num{0.48} \\
		specific heat capacity $\heatcapacity$ [\si{\joule \per \kilogram \per \kelvin}] & \num{3.69e3} \\
		heat transfer coefficient $\htc_\subcool$ [\si{\watt \per \kelvin \per \square\meter}] & \num{250} \\
		heat transfer coefficient $\htc_\subamb$ [\si{\watt \per \kelvin \per \square\meter}] & \num{0} \\
		coolant absorption factor $\coolingfactor$ [\si{1}] & \num{0.14} \\
		initial temperature $\temperature_0$ [\si{\celsius}] & \num{21.8} \\
		cooling temperature $\temperature_\subcool$ [\si{\celsius}] & \num{20} \\
		blood temperature $\temperature_\subblood$ [\si{\celsius}] & \num{21.8} \\
		ambient temperature $\temperature_\subamb$ [\si{\celsius}] & \num{21.8} \\
		\midrule
		optical (native) & \\
		\midrule
		absorption coefficient $\absorption{\subnative}$ [\si{1\per \meter}] & \num{50} \\
		scattering coefficient $\scattering{\subnative}$ [\si{1\per \meter}] & \num{8e3} \\
		anisotropy factor $\anisotropy{\subnative}$ [\si{1}] & \num{0.97} \\
		\midrule
		optical (coagulated) & \\
		\midrule
		absorption coefficient $\absorption{\subcoag}$ [\si{1 \per \meter}] & \num{60} \\
		scattering coefficient $\scattering{\subcoag}$ [\si{1 \per \meter}] & \num{3e4} \\
		anisotropy factor $\anisotropy{\subcoag}$ [\si{1}] & \num{0.95} \\
		\midrule
		tissue damage & \\
		\midrule
		frequency factor $\frequencyfactor$ [\si{1 \per \second}] & \num{3.1e98} \\
		activation energy $\activationenergy$ [\si{\joule \per \mole}] & \num{6.28e5}\\
		\bottomrule
	\end{tabular}
	\caption{Values for the various parameters of the liver tissue.}
	\label{table:parameters}
\end{table}

Finally, we consider the coagulation of the tissue molecules using Arrhenius' law (cf. \cite{fasano}) which reads
\begin{equation}
	\label{eq:arrhenius}
	\damage(t, x) = \int_{0}^{t} \frequencyfactor \exp\left( - \frac{\activationenergy}{\gasconstant \temperature(\theta, x)} \right) \dx{\theta},
\end{equation}
where $\frequencyfactor$ denotes the frequency factor [\si{1\per \second}], $\activationenergy$ is the activation energy [\si{\joule \per \mole}] and $\gasconstant$ the universal gas constant [\si{\joule \per \mole \per \kelvin}]. It is important to note that this quantity does not only depend on the temperature at the current time $t$, but on its entire history.
As tissue damage influences the optical parameters substantially, we model the effect of coagulation on them according to \cite{fasano, huebner} by
\begin{eqsystem*}
	\absorption{} &= \absorption{\subnative} + \Big(1 - \exp(-\damage)\Big) (\absorption{\subcoag} - \absorption{\subnative}), \\
	\scattering{} &= \scattering{\subnative} + \Big(1 - \exp(-\damage)\Big) (\scattering{\subcoag} - \scattering{\subnative}), \\
	\anisotropy{} &= \anisotropy{\subnative} + \Big(1 - \exp(-\damage)\Big) (\anisotropy{\subcoag} - \anisotropy{\subnative}),
\end{eqsystem*}
where the subscript $n$ denotes the native value of that parameter and the subscript $c$ the respective coagulated one. The parameter values used in this paper are shown in Table~\ref{table:parameters} (cf. \cite{huebner, puccini2003simulations, roggan1995optical, giering1995review, schwarzmaier1998treatment}). Note, that these parameters represent healthy tissue. It is, however, easily possible to account for tumorous tissue through local variations of the corresponding parameters.

\section{Parameter Identification}
\label{sec:parameter_identification}
We now formulate the parameter identification problem and relate it to an optimal control problem, which is subsequently treated with the help of the adjoint approach.

\subsection{Problem Formulation}

Given a temperature measurement $\temperature_\subdes$ at a certain time we want to identify (or reconstruct) the blood perfusion rate $\perfusion$ that induced the measured temperature distribution. Without loss of generality we assume that the time of the measurement coincides with the time horizon $\timehorizon$. Later on, we choose a time horizon that is much smaller than the end time of the therapy $\timehorizon_\subend$ to identify the perfusion rate at the beginning of the treatment. We assume that our state equations \eqref{eq:bioheat} and \eqref{eq:radiation} admit a unique solution (cf. \cite{tse}). We define the cost functional $J(\temperature, \perfusion)$ as
\begin{equation}
	\label{eq:def_cost_function}
	J(\temperature, \perfusion) := \frac{1}{2} \int_{\Omega} \abs{\temperature(\timehorizon, \cdot) - \temperature_\subdes}^2 \dx{x} + \frac{\lambda}{2} \int_{\Omega} \abs{\perfusion}^2 \dx{x},
\end{equation}
which is then used to model the parameter identification problem described above with the following minimization problem:
\begin{equation}
	\label{eq:optimal_control_problem}
	\min_{\perfusion \in \mathcal{A}, \temperature}\ J(\temperature, \perfusion) \quad \text{ s.t. } \eqref{eq:bioheat} \text{ and } \eqref{eq:radiation}.
\end{equation}
Here, the perfusion rate $\perfusion$ plays the role of a control that is used to \qe{steer} the state, i.e., the simulated temperature, to the desired state, i.e., the temperature measurement. Further, $\mathcal{A} \subset L^2(\Omega)$ denotes the set of admissible perfusion rates which is used in order to model so-called control constraints, e.g., only nonnegative perfusion rates are physically meaningful. The first term (also known as observation term) tries to minimize the difference of $\temperature(\tau, \cdot)$ and $\temperature_\subdes$. This means that we try to compute a perfusion rate such that the resulting temperature distribution at time $\timehorizon$ is close to the measured temperature. The second term is a so-called Tikhonov regularization with regularization parameter $\lambda \geq 0$, which is used to stabilize the (possibly) ill-posed problem (cf. \cite{rieder, Poerner2018}).

We reformulate the optimization problem \eqref{eq:optimal_control_problem} equivalently thanks to our assumption that the state equations admit a unique solution. To do this, we denote by $\temperature[\perfusion]$ the solution of \eqref{eq:bioheat} and \eqref{eq:radiation} with blood perfusion rate $\perfusion$. We introduce the reduced cost functional $\hat{J}$ by
\begin{equation}
	\label{eq:def_reduced_cost_functional}
	\hat{J}(\perfusion) := J(\temperature[\perfusion], \perfusion),
\end{equation}
and see that \eqref{eq:optimal_control_problem} is equivalent to the reduced problem
\begin{equation}
	\label{eq:reduced_optimal_control_problem}
	\min_{\perfusion\in \mathcal{A}}\ \hat{J}(\perfusion) = \frac{1}{2} \int_{\Omega} \Big\lvert\temperature[\perfusion](\timehorizon, \cdot) - \temperature_\subdes\Big\rvert^2 \dx{x} + \frac{\lambda}{2} \int_{\Omega} \abs{\perfusion}^2 \dx{x},
\end{equation}
where the PDE constraint is formally eliminated. To solve this minimization problem, we apply techniques from PDE-constrained optimization. In particular, we compute the gradient of the reduced cost functional $\hat{J}$ which is then used to solve problem \eqref{eq:reduced_optimal_control_problem} numerically with a gradient descent or a quasi-Newton method (cf. Section~\ref{sec:numerics}). For a detailed introduction to optimization problems constrained by PDEs and their (numerical) solution we refer, e.g., to \cite{troeltzsch, hinze_pinnau_ulbrich2, borzi_schulz}.

\subsection{Adjoint-Based Identification}

To compute the gradient of $\hat{J}$, we use the formal Lagrange method of \cite[Chapter 2.10]{troeltzsch}. For this purpose, we set up a Lagrangian $\mathcal{L}(\perfusion, \temperature, \radiation, \adtemp, \adrad)$, where $\adtemp$ and $\adrad$ are used as Lagrange multipliers for the PDE constraints \eqref{eq:bioheat} and \eqref{eq:radiation}. Then, the first order optimality conditions for a stationary point of the Lagrangian (and, therefore, for a minimizer of \eqref{eq:optimal_control_problem}) are given by the system
\begin{alignat}{2}
	\label{eq:kkt_state_bioheat}
	\frac{\partial \mathcal{L}}{\partial \adtemp}(\perfusion, \temperature, \radiation, \adtemp, \adrad)[\hat{\adtemp}] =& 0 \quad &&\text{ for all } \hat{\adtemp}, \\
	 \label{eq:kkt_state_radiation}
	\frac{\partial \mathcal{L}}{\partial \adrad}(\perfusion, \temperature, \radiation, \adtemp, \adrad)[\hat{\adrad}] =& 0 \quad &&\text{ for all } \hat{\adrad}, \\
	\label{eq:kkt_adjoint_bioheat}
	\frac{\partial \mathcal{L}}{\partial \temperature}(\perfusion, \temperature, \radiation, \adtemp, \adrad)[\hat{\temperature}] =& 0 \quad &&\text{ for all } \hat{\temperature}, \\
	\label{eq:kkt_adjoint_radiation}
	\frac{\partial \mathcal{L}}{\partial \radiation}(\perfusion, \temperature, \radiation, \adtemp, \adrad)[\hat{\radiation}] =& 0 \quad &&\text{ for all } \hat{\radiation},\quad \text{ and } \\
	\label{eq:kkt_optimality_condition}
	\frac{\partial \mathcal{L}}{\partial \perfusion}(\perfusion, \temperature, \radiation, \adtemp, \adrad)[\hat{\perfusion} - \perfusion] 
	\geq& 0 \quad &&\text{ for all } \hat{\perfusion} \in \mathcal{A}.
\end{alignat}

For our problem \eqref{eq:optimal_control_problem}, this Lagrangian is given by
\begin{eqsystem}[eq:lagrangian]
	\mathcal{L}(\perfusion, \temperature, \radiation, \adtemp, \adrad) = &\frac{1}{2} \int_{\Omega} \abs{\temperature(\timehorizon, \cdot) - \temperature_\subdes}^2 \dx{x} + \frac{\lambda}{2} \int_{\Omega} \abs{\perfusion}^2 \dx{x} \\
	- &\int_{0}^{\timehorizon} \int_{\Omega} \density \heatcapacity \dot{\temperature} \adtemp_1 \dx{x} \dx{t} + \int_{0}^{\timehorizon} \int_{\Omega} \nabla \cdot (\conductivity \nabla \temperature)\ \adtemp_1 \dx{x} \dx{t} \\
	- &\int_{0}^{\timehorizon} \int_{\Omega} \perfusion (\temperature - \temperature_\subblood) \adtemp_1 \dx{x} \dx{t} + \int_{0}^{\timehorizon} \int_{\Omega} \absorption{} \radiation \adtemp_1 \dx{x} \dx{t} \\
	- &\int_{0}^{\timehorizon} \int_{\Gamma} \conductivity \normal \cdot \nabla \temperature\ \adtemp_2 \dx{s} \dx{t} - \int_{0}^{\timehorizon} \int_{\Gamma} \htc (\temperature - \temperature_\subext)\ \adtemp_2 \dx{s} \dx{t} \\
	- &\int_{\Omega} (\temperature(0, \cdot) - \temperature_0) \adtemp_3 \dx{x} + \int_{0}^{\timehorizon} \int_{\Omega} \nabla \cdot (D \nabla \radiation) \adrad_1 \dx{x} \dx{t} \\
	- &\int_{0}^{\timehorizon} \int_{\Omega} \absorption{} \radiation \adrad_1 \dx{x} \dx{t} - \int_{0}^{\timehorizon} \int_{\Gamma_\subrad} D \normal \cdot \nabla \radiation\ \adrad_2 \dx{s} \dx{t} \\
	+ &\int_{0}^{\timehorizon} \int_{\Gamma_\subrad} \frac{\laserpower}{\abs{\Gamma_\subrad}} \adrad_2 \dx{s} \dx{t} - \int_{0}^{\timehorizon} \int_{\Gamma_\subcool} D \normal \cdot \nabla \radiation\ \adrad_3 \dx{s} \dx{t} \\
	- &\int_{0}^{\timehorizon} \int_{\Gamma_\subamb} D \normal \cdot \nabla \radiation\ \adrad_4 \dx{s} \dx{t} - \int_{0}^{\timehorizon} \int_{\Gamma_\subamb} \frac{1}{2} \radiation \adrad_4 \dx{s} \dx{t},
\end{eqsystem}
where $\adtemp = [\adtemp_1, \adtemp_2, \adtemp_3]$ and $\adrad = [\adrad_1, \adrad_2, \adrad_3, \adrad_4]$. As the Lagrangian is linear in $\adtemp$ and $\adrad$, we see that equations \eqref{eq:kkt_state_bioheat} and \eqref{eq:kkt_state_radiation} of the optimality system are equivalent to the state equations \eqref{eq:bioheat} and \eqref{eq:radiation} that constrain the optimization problem \eqref{eq:optimal_control_problem}.

Lengthy calculations show that \eqref{eq:kkt_adjoint_bioheat} and \eqref{eq:kkt_adjoint_radiation} give rise to the conditions
\begin{align*}
	&\adrad_2 = \adrad_1 \quad \text{ on } (0, \timehorizon) \times \Gamma_\subrad, \qquad \adrad_3 = \adrad_1 \quad \text{ on } (0, \timehorizon) \times \Gamma_\subcool, \\
	&\adrad_4 = \adrad_1 \quad \text{ on } (0, \timehorizon) \times \Gamma_\subamb, \qquad \adtemp_2 = \adtemp_1 \quad \text{ on } (0, \timehorizon) \times \Gamma, \quad \text{ and } \\
	&\adtemp_3 = \density \heatcapacity \adtemp_1(0),
\end{align*}
such that we only have to consider the multipliers $\adtemp_1$ and $\adrad_1$. Therefore, in the following we drop the index and only write $\adtemp = \adtemp_1$ as well as $\adrad = \adrad_1$. Furthermore, with the above conditions, \eqref{eq:kkt_adjoint_bioheat} and \eqref{eq:kkt_adjoint_radiation} are equivalent to the following system of adjoint equations:
\begin{eqsystem}[eq:adjoint_bioheat]
	- \density \heatcapacity \dot{\adtemp} - \nabla \cdot (\conductivity \nabla \adtemp) + \perfusion \adtemp + F_1 + F_2 &= 0 \quad &&\text{ in } (0, \timehorizon) \times \Omega, \\
	\conductivity \partial_\normal \adtemp + \htc \adtemp &= 0 \quad &&\text{ on } (0, \timehorizon) \times \Gamma, \\
	\adtemp(\timehorizon, \cdot) - \frac{\lambda}{\density \heatcapacity} \left( \temperature(\timehorizon, \cdot) - \temperature_\subdes \right) &= 0 \quad &&\text{ in } \Omega,
\end{eqsystem}
where 
\begin{align*}
	F_1 =& \frequencyfactor \exp\left( -\frac{\activationenergy}{\gasconstant \temperature}  \right) \frac{\activationenergy}{\gasconstant \temperature^2} \int_{t}^{\timehorizon} \frac{\partial \absorption{}}{\partial \damage} \radiation \left( \adrad - \adtemp \right) \dx{\theta}, \\
	F_2 =& \frequencyfactor \exp\left( -\frac{\activationenergy}{\gasconstant \temperature}  \right) \frac{\activationenergy}{\gasconstant \temperature^2} \int_{t}^{\timehorizon} \frac{\partial \diffusivity}{\partial \damage} \nabla \radiation \cdot \nabla \adrad \dx{\theta},
\end{align*}
and $\adtemp$ solves
\begin{eqsystem}[eq:adjoint_radiation]
	-\nabla \cdot (D \nabla \adrad) + \absorption{} \adrad &= \absorption{} \adtemp \quad &&\text{ in } (0, \timehorizon) \times \Omega, \\
	D \partial_\normal \adrad &= 0 \quad &&\text{ on } (0, \timehorizon) \times \left( \Gamma_\subcool \cup \Gamma_\subrad \right), \\
	D \partial_\normal \adrad + \frac{1}{2} \adrad &= 0 \quad &&\text{ on } (0, \timehorizon) \times \Gamma_\subamb.
\end{eqsystem}
Moreover, we remark that, as usual, the adjoint (bio-)heat equation is an equation that is posed backward in time. For the analysis and numerical solution of this equation one can introduce the time shift $\Theta = \tau - t$ and then solve an equation that evolves forward in $\Theta$ (cf. \cite{troeltzsch, hinze_pinnau_ulbrich2}). 

Finally, the optimality condition \eqref{eq:kkt_optimality_condition} is equivalent to the following variational inequality 
\begin{equation}
	\label{eq:optimality_condition}
	\hat{J}'(\perfusion)\left[\hat{\perfusion} - \perfusion\right] = \left( \hat{J}'(\perfusion), \hat{\perfusion} - \perfusion \right)_{L^2(\Omega)} \geq 0 \quad \text{ for all } \hat{\perfusion} \in \mathcal{A},
\end{equation}
where $\left( u, v \right)_H$ denotes the inner product of $u$ and $v$ in some Hilbert space $H$. Here, the gradient of the reduced cost functional is given by
\begin{equation}
	\label{eq:gradient}
	\hat{J}'(\perfusion) = \lambda \perfusion + \int_{0}^{\timehorizon} \left( \temperature_\subblood - \temperature \right) \adtemp \dx{t}.
\end{equation}
This relation is used in Section~\ref{sec:numerics} for the numerical solution of the parameter identification problem \eqref{eq:optimal_control_problem}. Our results are summarized in Proposition~\ref{prop:gradient}.

\begin{Proposition}
	\label{prop:gradient}
	Let $\mathcal{A}$ be a convex subset of $L^2(\Omega)$. The first order necessary conditions for $\perfusion^*$ being a minimizer of \eqref{eq:reduced_optimal_control_problem} are given by the state equations \eqref{eq:bioheat} and \eqref{eq:radiation}, the adjoint equations \eqref{eq:adjoint_bioheat} and \eqref{eq:adjoint_radiation}, as well as the variational inequality \eqref{eq:optimality_condition}. The gradient $\hat{J}'(\perfusion)$ of the reduced cost functional is given by \eqref{eq:gradient}.	
\end{Proposition}

\section{Numerical Methods}
\label{sec:numerics}

In the following we discuss the numerical methods used to solve the parameter identification problem. First, we describe the numerical solution techniques for the state and adjoint equations and then we elaborate the algorithms used for solving the optimization problem.

\subsection{Solution of the PDEs}

We solve all PDEs, i.e., the state and adjoint equations, with the finite element method. For this purpose, we triangulate our domain $\Omega$ with the help of GMSH, version 2.11.0 \cite{gmsh}. The assembly and solution of the linear systems is done with FEniCS, version 2018.1 \cite{fenics, fenics_book} and PETSc, version 3.10.5 \cite{petsc-user-ref}, respectively. To solve the time-dependent PDEs \eqref{eq:bioheat} and \eqref{eq:adjoint_bioheat} we first discretize them in time with the implicit Euler method. Further, all PDEs are discretized in space with the help of linear Lagrange elements. The resulting sequences of linear systems corresponding to \eqref{eq:bioheat}, \eqref{eq:radiation}, and \eqref{eq:adjoint_radiation} are then solved with the conjugate gradient method and an incomplete Cholesky factorization as preconditioner. For the solution of the sequence of linear systems arising from \eqref{eq:adjoint_bioheat} we use the \texttt{MINRES} algorithm and the Jacobi method as preconditioner, as the corresponding matrices are symmetric, but not necessarily positive definite.

\subsection{Solution of the Optimization Problem}

Let us now turn our attention to the numerical solution of the optimal control problem \eqref{eq:optimal_control_problem} which we solve by the means of a projected quasi-Newton method described in \cite{kelley}. The idea of the method is the following: Assume that we computed the $k$-th iterate $\perfusion_k$. To compute the gradient of the reduced cost functional, we first solve the state equations \eqref{eq:bioheat} and \eqref{eq:radiation} and then the adjoint equations \eqref{eq:adjoint_bioheat} and \eqref{eq:adjoint_radiation}. Subsequently, we compute $g_k = \hat{J}'(\perfusion_k)$ with the relation \eqref{eq:gradient} and, with this, the search direction $d_k$ is computed by a L-BFGS update of the form
\begin{equation}
	\label{eq:quasi_newton}
	d_k = -H_k^{-1} g_k,
\end{equation}
where $H_k$ denotes the L-BFGS approximation of the (reduced) Hessian of $\hat{J}$ at $\perfusion_k$. The application of the inverse of $H_k$ is efficiently performed via the well-known BFGS double loop \cite{kelley, nocedal_wright}. Next, we perform a line search along the projected path given by $\mathcal{P}(\perfusion_k + \alpha d_k)$, where $\mathcal{P}$ denotes the projection onto $\mathcal{A}$, to find a suitable step size $\alpha_k$. This is done by the following Armijo rule (cf. \cite{bertsekas, calamai_more}): Define $\perfusion_k(\alpha) = \mathcal{P}(\perfusion_k + \alpha d_k)$. Then, the step size $\alpha_k$ is of the form $\alpha_k = \beta^{m_k} \gamma$, where $\beta \in (0, 1)$ and $\gamma >0$ and $m_k$ is the smallest integer satisfying 
\begin{equation}
	\label{eq:armijo}
	\hat{J}(\perfusion_k(\alpha_k)) \leq \hat{J}(\perfusion_k) + \gamma \left( \hat{J}'(\perfusion_k), \perfusion_k(\alpha_k) - \perfusion_k \right)_{L^2(\Omega)}.
\end{equation}
For all of our numerical results we choose $\beta = \nicefrac{1}{2}$ and $\gamma = \num{1e-4}$ (cf. \cite{nocedal_wright}). Finally, we update the iterate by $\perfusion_{k+1} = \mathcal{P}(\perfusion_k + \alpha_k d_k)$. For the stopping criterion we define the stationary measure as 
\begin{equation*}
	\Sigma (\perfusion) := \norm{\perfusion - \mathcal{P}(\perfusion - \hat{J}'(\perfusion))}{L^2(\Omega)},
\end{equation*}
 and terminate the method once the relative stopping criterion
\begin{equation}
	\label{eq:stopping_criterion}
	\Sigma(\perfusion_k) \leq \texttt{tol}\ \Sigma(\perfusion_0)
\end{equation}
is satisfied, where \texttt{tol} is a user-defined tolerance (cf. \cite{kelley}). This procedure is summarized in Algorithm~\ref{algo:descent}. A detailed description of the optimization methods can be found in, e.g., \cite{kelley, nocedal_wright}, and their application to PDE-constrained optimization problems is covered in, e.g., \cite{hinze_pinnau_ulbrich2, troeltzsch, borzi_schulz}.

\begin{algorithm}
	\KwIn{initial perfusion rate $\perfusion_0$, tolerance $\texttt{tol} \in (0, 1)$}
	
	\For{$k = 0,1,2, \dots$}{
		\If{Stopping criterion \eqref{eq:stopping_criterion} is satisfied}{Stop with minimizer $\perfusion_k$}
		Solve the state equations \eqref{eq:bioheat} and \eqref{eq:radiation} \\
		Solve the adjoint equations \eqref{eq:adjoint_bioheat} and \eqref{eq:adjoint_radiation} \\
		Compute the gradient $g_k = \hat{J}'(\perfusion_k)$ via \eqref{eq:gradient}\\
		Compute the search direction $d_k$ as $H_kd_k = - \hat{J}'(\perfusion_k)$ \label{algo:computation} \\
		Compute a feasible stepsize $\alpha_k$ with the Armijo rule \eqref{eq:armijo}\\
		Update the perfusion rate: $\perfusion_{k+1} = \mathcal{P}(\perfusion_k + \alpha_k d_k) $
	}
	\caption{Projected quasi-Newton algorithm.}
	\label{algo:descent}
\end{algorithm}

For the computation of the search direction with the quasi-Newton method in step \ref{algo:computation} of Algorithm~\ref{algo:descent} (cf. \eqref{eq:quasi_newton}) we use the algorithm given in \cite[Chapter 5.5.3]{kelley}. This corresponds to a projected BFGS method that approximates the reduced Hessian of the problem. We implemented this by the means of a limited memory BFGS update that only uses the information of the last $l \in \mathbb{N}$ steps. In particular, we get a complete projected BFGS method in case $l = \infty$. On the other hand, for $l=0$ we choose $H_k = I$, where $I$ denotes the identity, and the whole method reduces to the projected gradient descent method that is described, e.g., in \cite{kelley, hinze_pinnau_ulbrich2, troeltzsch}. Therefore, when we speak of using a (projected) gradient descent method we refer to the case $l=0$ and when we talk about using a (projected) L-BFGS method we refer to the case $l > 0$. These two methods are compared in Section~\ref{sec:results}. Finally, note that in case the curvature condition for the BFGS method is not satisfied we re-initialize the method with the identity, i.e., we perform a gradient descent step (see \cite[Chapter 4.2.2]{kelley} for details).

\subsection{Multiple Measurements}

The ideas and methods described before can be generalized to the case where multiple measurements are taken during the therapy. To do so, assume that $m$ measurements $\temperature_\subdes^{(1)}, \dots, \temperature_\subdes^{(m)}$ are taken at times $\timehorizon_1 < \timehorizon_2 < \dots < \timehorizon_m$ and that we have an initial guess $\perfusion_0$ for the perfusion rate. As before, all measurements should be taken before the end of the therapy, so $\timehorizon_m < \timehorizon_\subend$ and, additionally, we define $\timehorizon_0 = 0$.

A first approach for solving this problem would be to use Algorithm~\ref{algo:descent} on each interval $(\timehorizon_k, \timehorizon_{k+1})$ separately, where the initial temperature is given by the measurement. However, this has the important drawback, that we would also need measurements of tissue damage at $\timehorizon_k$ which we cannot compute accurately from the thermometry data. Therefore, we consider using the previously simulated temperature distribution and tissue damage as initial conditions for the subsequent identification interval. 

In particular, our method proceeds as follows: On the first interval, i.e., $(0, \timehorizon_1)$ we use Algorithm~\ref{algo:descent} in order to identify the blood perfusion rate $\perfusion^{(1)}$ and the resulting temperature distribution $\temperature^{(1)} = \temperature[\perfusion^{(1)}]$ as well as damage function $\omega^{(1)}$. These are then used as initial conditions for the state equations for the identification in the subsequent interval $(\timehorizon_1, \timehorizon_2)$. Furthermore, we also use the perfusion rate $\perfusion^{(1)}$ as initial guess for the parameter identification algorithm on $(\timehorizon_1, \timehorizon_2)$. These steps are then repeated until the last identification is carried out. Therefore, the blood perfusion rate we compute with this approach is constant on each interval $(\timehorizon_k, \timehorizon_{k+1})$, and the simulated temperature is the one corresponding to this piecewise constant perfusion rate. As we want to predict the temperature in the time interval $(\timehorizon_m, \timehorizon_\subend)$, we choose the last computed perfusion rate $\perfusion^{(m)}$ as the perfusion rate on this interval. We do so as the effect of the perfusion rate depends on the magnitude of the difference $\temperature - \temperature_\subblood$ which increases over time. Therefore, the last identified perfusion rate should also be the most accurate one. The numerical experiments described in the following section confirm that this approach works well.

\section{Numerical Results}
\label{sec:results}

\begin{figure}[!t]
	\centering
	\includegraphics[width=\textwidth]{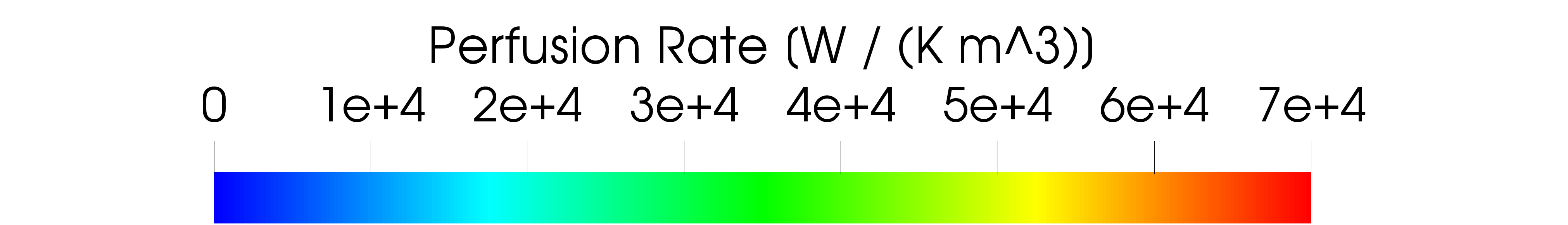}
	\begin{subfigure}[t]{0.45\textwidth}
		\includegraphics[width=\textwidth, trim=17.5cm 2.5cm 17.5cm 2.5cm, clip]{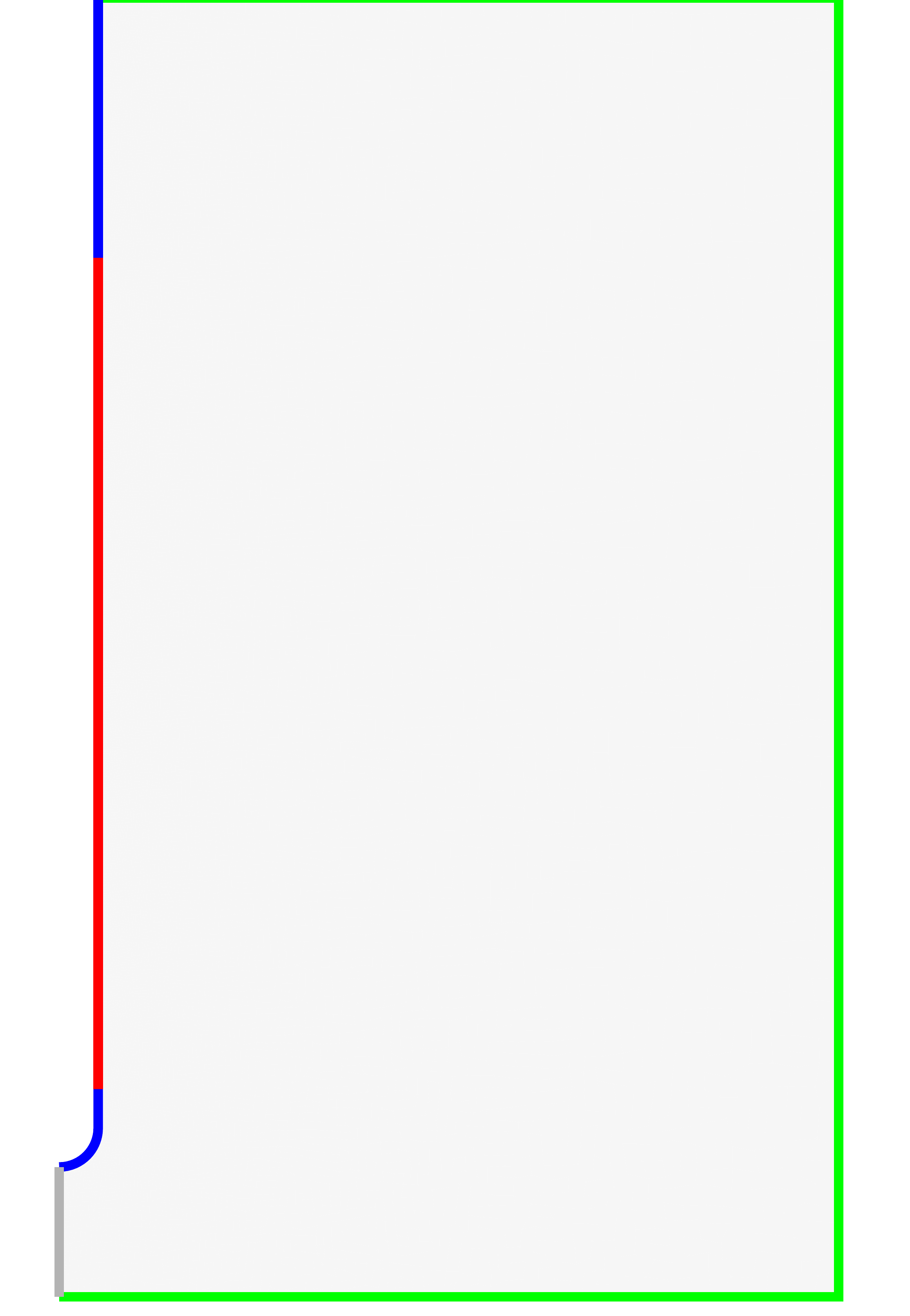}
		\caption{Axisymmetric geometry $\Omega$ with radiating boundary $\Gamma_\subrad$ (red), cooling boundary $\Gamma_\subcool$ (blue), ambient boundary $\Gamma_\subamb$ (green) and symmetry boundary (gray).}
		\label{fig:geometry}
	\end{subfigure}
	\hfill
	\begin{subfigure}[t]{0.45\textwidth}
		\includegraphics[width=0.9\textwidth]{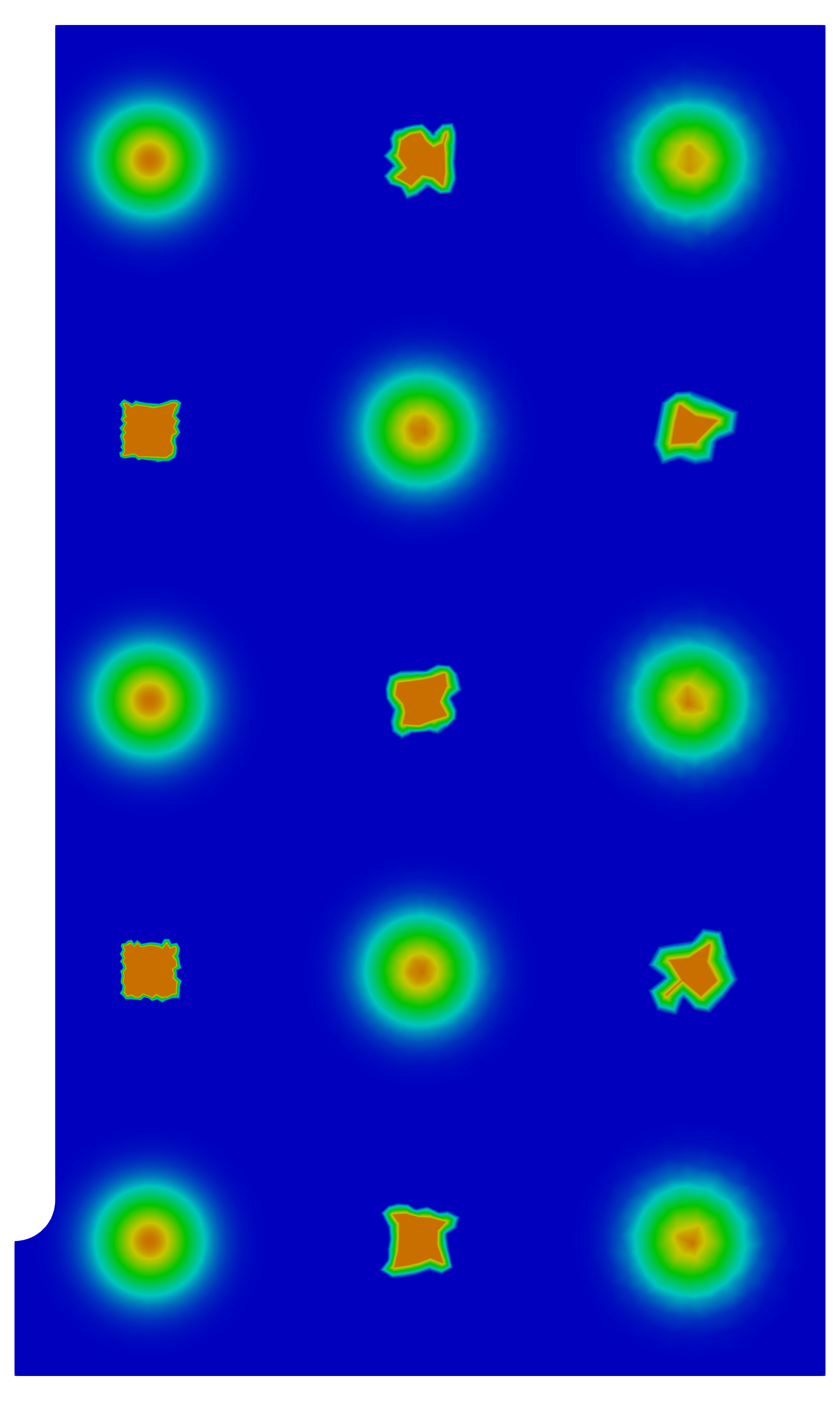}
		\caption{Prescribed perfusion rate.}
		\label{fig:prescribed}
	\end{subfigure}
	\caption{Problem setting.}
\end{figure}

We now apply the previously introduced techniques for identifying the perfusion rate to a model problem: We generate an artificial temperature measurement $\temperature_\subdes$ from a synthetic perfusion rate $\perfusion_\subdes$ as solution of the state equations. Furthermore, for all of our experiments we choose the set of admissible perfusion rates $\mathcal{A}$ as
\begin{equation*}
	\mathcal{A} := \Set{\xi \in L^2(\Omega) |  \xi \geq 0 \text{ a.e. in } \Omega }.
\end{equation*}
As stated previously, this is sensible as all physically meaningful perfusion rates are nonnegative. 

To demonstrate the behavior and capabilities of the parameter identification algorithm of Section~\ref{sec:numerics}, we first assume that there is no noise present in the measurement, i.e., we perform the identification with an \qe{exact} measurement. Afterwards, we investigate its performance in the presence of noisy measurements. For simplicity, consider axisymmetric problems such that we can perform the numerical studies in 2D. The axisymmetric geometry is shown in Figure~\ref{fig:geometry}. The parameters used for our problem are taken from \cite{huebner} and are depicted in Table~\ref{table:parameters}.
Note that these parameters represent ex-vivo porcine tissue. However, the values are close to the ones of human tissue \cite{puccini2003simulations} and the results can be transferred to the in-vivo scenario.

\subsection{Noiseless Model Problem}
\label{sec:numerics_noiseless}

For this problem, we only compute the perfusion rate after $\timehorizon_1 = \num{60}\ \si{\second}$ and simulate the rest of the treatment with the perfusion rate computed in this first identification step. For the optimization algorithms we choose the initial guess for the perfusion rate as $\perfusion = 0$ and the relative stopping tolerance is set to $\texttt{tol} = \num{1e-3}$. Furthermore, we set the regularization parameter as $\lambda = 0$, i.e., we do not use a Tikhonov regularization for this case. We stop the parameter identification after \num{20} iterations in case the stopping criterion is not satisfied. The synthetic perfusion rate $\perfusion_\subdes$ is depicted in Figure~\ref{fig:prescribed}. The maximum perfusion rate is chosen to be $\perfusion_\text{max} = \num{6e4}~\si{\watt \per \kelvin \per \cubic \meter}$ in accordance to \cite{perfusion_parameter}. Outside the blood vessels we set the perfusion rate to \num{0}. For the numerical analysis of the parameter identification we consider two different types of blood vessels: First, \qe{smooth vessels} that are modeled via two-dimensional Gaussian kernels with maximum height $\perfusion_\text{max}$, and second, \qe{square vessels} that have a constant perfusion rate $\perfusion_\text{max}$. Note that for the L-BFGS algorithm we use (at most) the information of the last \num{5} iterations for the update of the approximate Hessian.

\begin{figure}[!b]
	\centering
	\includegraphics[width=\textwidth]{legend}
	\begin{subfigure}[t]{0.475\textwidth}
		\includegraphics[width=\textwidth]{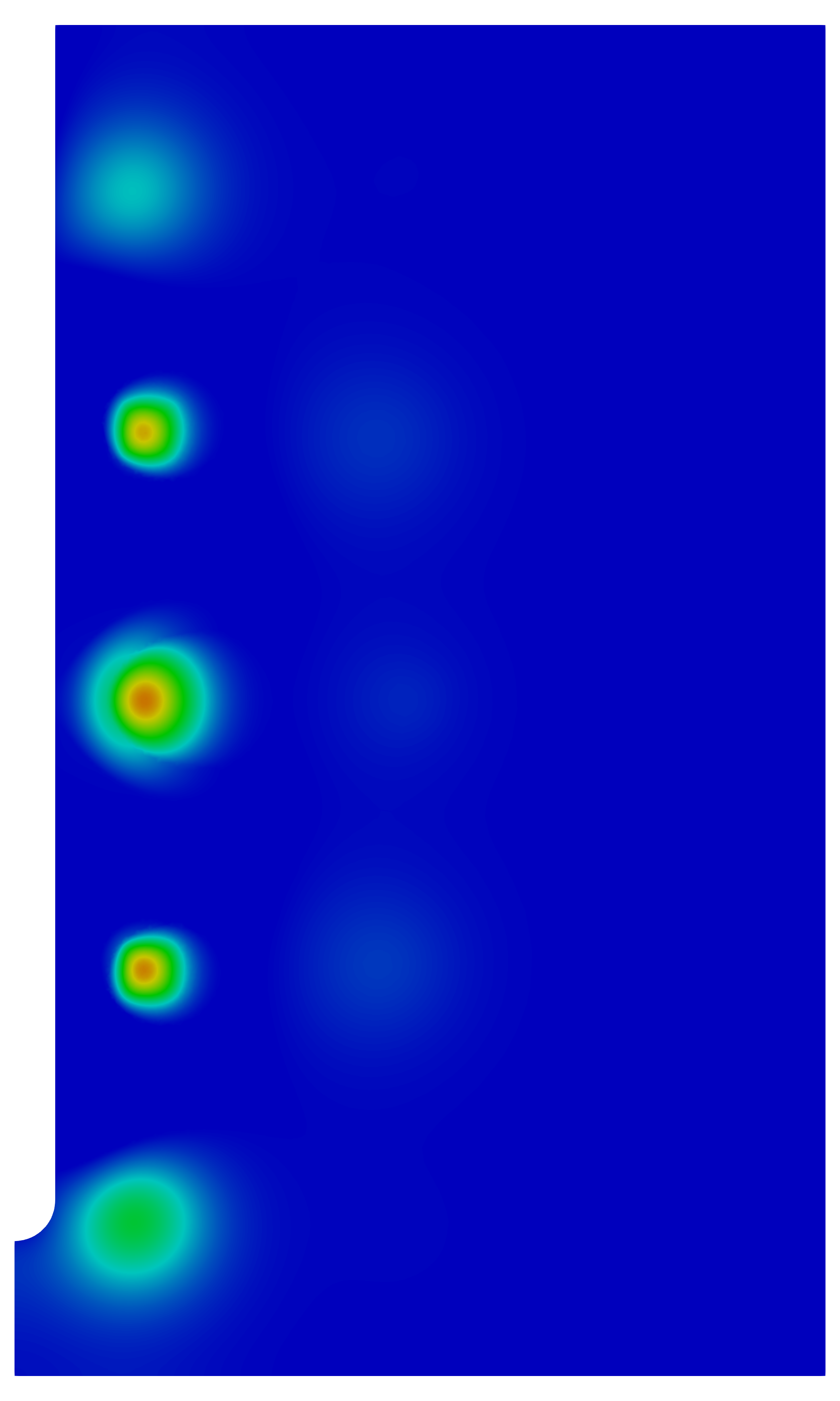}
		\caption{Projected gradient descent method.}
		\label{fig:identified_gradient}
	\end{subfigure}
	\hfil
	\begin{subfigure}[t]{0.475\textwidth}
		\includegraphics[width=\textwidth]{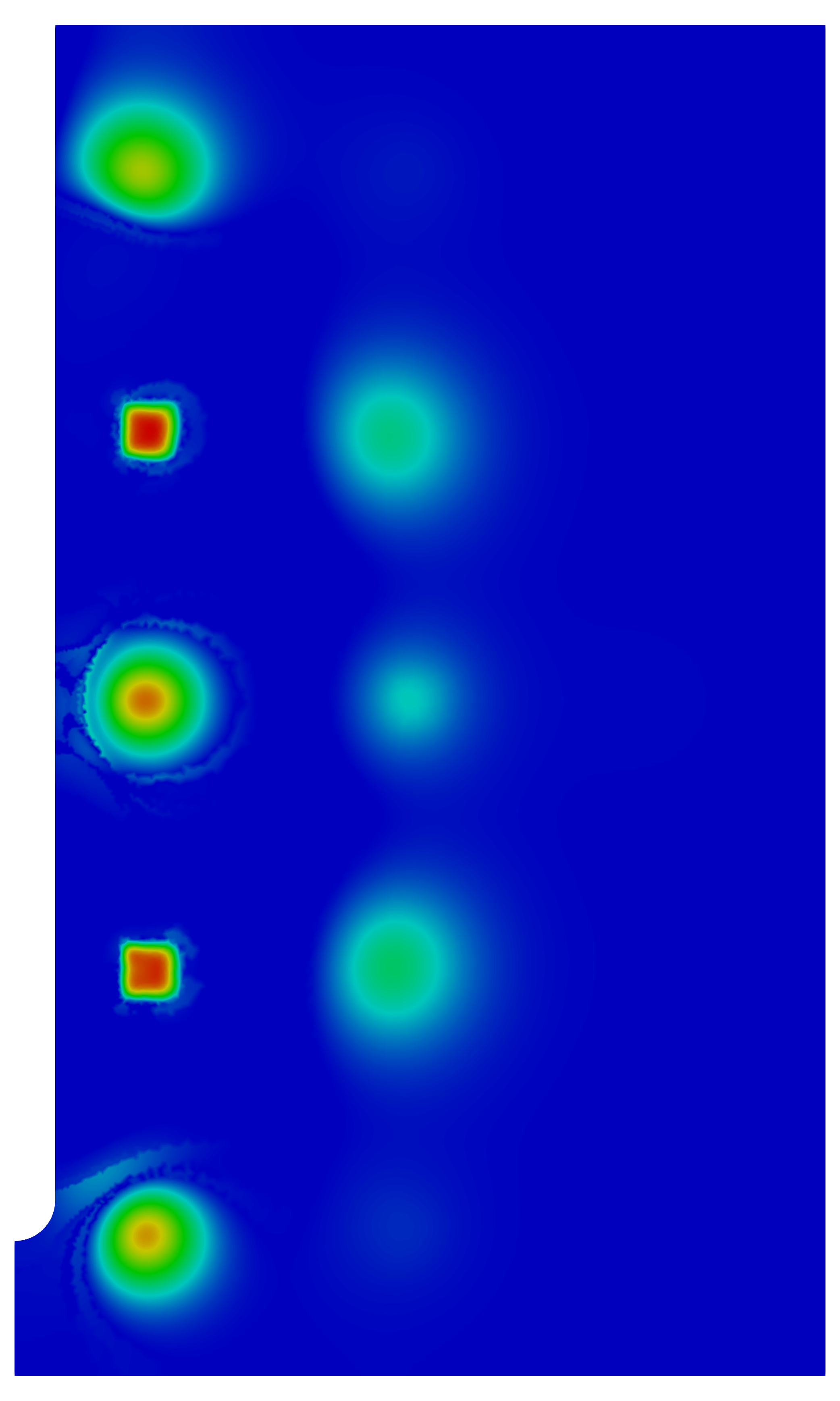}
		\caption{Projected L-BFGS method.}
		\label{fig:identified_bfgs}
	\end{subfigure}
	\caption{Identified perfusion rates.}
	\label{fig:identified_noiseless}
\end{figure}

\begin{figure}[!t]
	\begin{subfigure}{0.475\textwidth}
		\includegraphics[width=\textwidth]{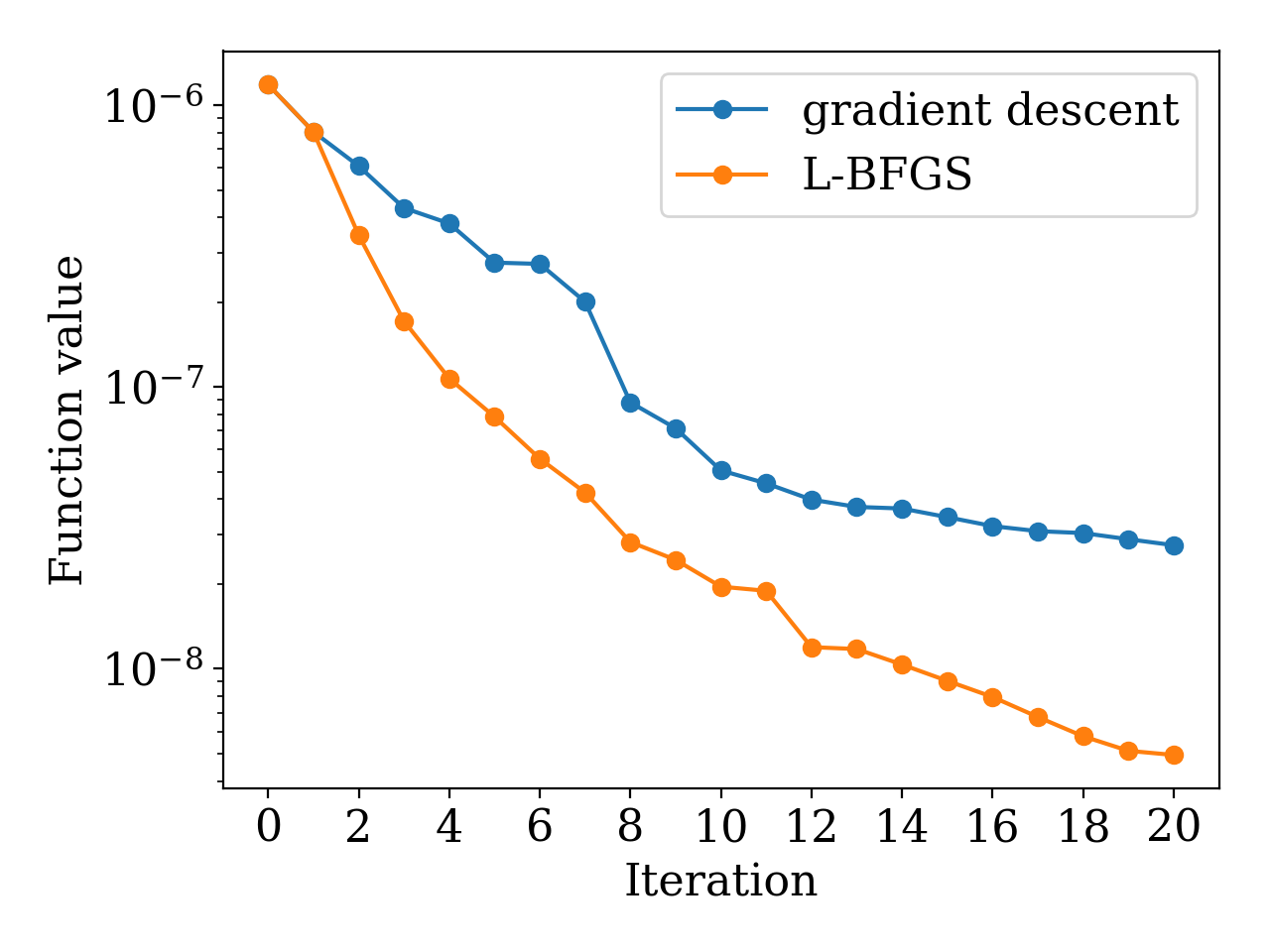}
		\caption{Evolution of $\hat{J}(\perfusion_k)$.}
		\label{fig:history_function_noiseless}
	\end{subfigure}
	\hfil
	\begin{subfigure}{0.475\textwidth}
		\includegraphics[width=\textwidth]{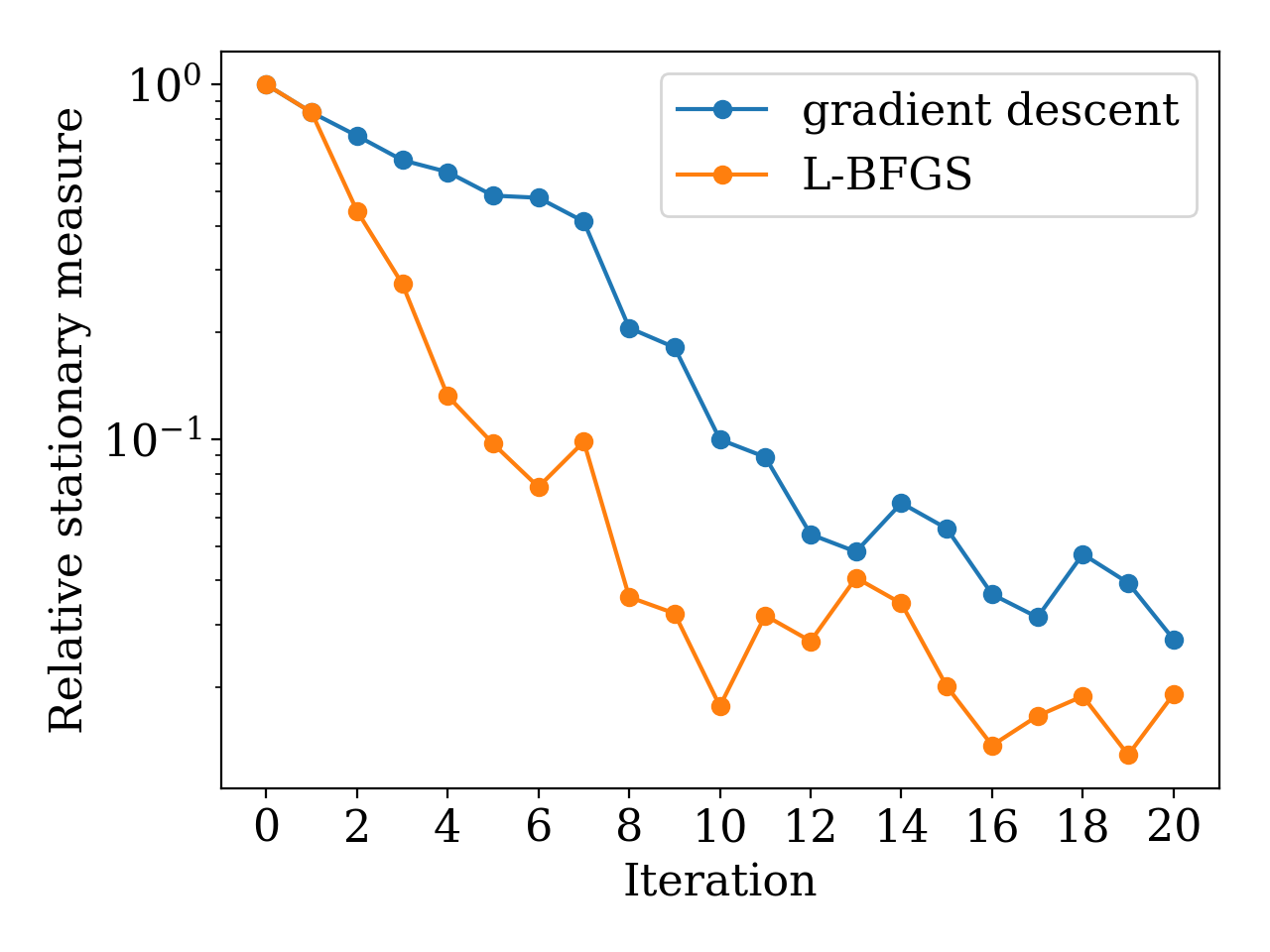}
		\caption{Evolution of $\nicefrac{\Sigma(\perfusion_k)}{\Sigma(\perfusion_0)}$.}
		\label{fig:history_measure_noiseless}
	\end{subfigure}
	\caption{Convergence history of the optimization methods.}
	\label{fig:history_noiseless}
\end{figure}

\begin{table}[!b]
	\centering
	\begin{tabular}{l r r }
		\toprule
		& $L^\infty$-error & $L^2$-error  \\
		\midrule
		gradient descent & & \\
		\midrule
		temperature $\temperature$ & \num{8.684} (\num{2.111} \%) & \num{0.916} (\num{0.731} \%) \\
		radiative energy $\radiation$ & \num{9931} (\num{2.738} \%) &  \num{27.711} (\num{0.456} \%) \\
		tissue damage $\dmgfun$ & \num{0.863} (\num{86.26} \%) & \num{0.017} (\num{19.05} \%) \\
		\midrule
		L-BFGS & & \\
		\midrule
		temperature $\temperature$ & \num{3.789} (\num{0.921} \%) & \num{0.413} (\num{0.33} \%) \\
		radiative energy $\radiation$ & \num{7002} (\num{1.93} \%) & \num{9.092} (\num{0.15} \%) \\
		tissue damage $\dmgfun$ & \num{0.253} (\num{25.295} \%) & \num{0.002} (\num{1.8} \%) \\
		\midrule
		$\perfusion = 0$ & & \\
		\midrule
		temperature $\temperature$ & \num{44.835} (\num{10.9} \%) & \num{2.629} (\num{2.098} \%) \\
		radiative energy $\radiation$ & \num{2.103e5} (\num{57.985} \%) & \num{372.502} (\num{6.129} \%) \\
		tissue damage $\dmgfun$ & \num{1} (\num{100} \%) & \num{0.063} (\num{69.422} \%) \\
		\bottomrule
	\end{tabular}
	\caption{Comparison between the simulated and measured temperature for a single measurement.}
	\label{table:comparison_noiseless}
\end{table}

\subsection*{A Single Measurement}
The identified perfusion rates are depicted in Figure~\ref{fig:identified_noiseless}, where the result of the gradient descent method is shown in Figure~\ref{fig:identified_gradient} and the perfusion rate computed with the L-BFGS method can be seen in Figure~\ref{fig:identified_bfgs}. The convergence history of both optimization algorithms is shown in Figure~\ref{fig:history_noiseless}, where the function values are given in Figure~\ref{fig:history_function_noiseless} and the relative stationary measure is depicted in Figure~\ref{fig:history_measure_noiseless}. Finally, we also compare the simulated temperature, radiative energy and tissue damage to our artificial measurements. The absolute and relative errors in $L^\infty(0,\timehorizon_\subend; L^\infty(\Omega))$ and $L^2(0,\timehorizon_\subend; L^2(\Omega))$ norms can be seen in Table~\ref{table:comparison_noiseless} for both the gradient descent method and the L-BFGS method. Note that we define the tissue damage as $\dmgfun = 1 - \exp(-\omega)$, i.e., the measure that indicates whether the tissue is in its native ($\dmgfun = 0$) or coagulated ($\dmgfun = 1$) state as the parameter $\damage$ enters our model only through $\dmgfun$.

\begin{figure}[!t]
	\begin{subfigure}{0.475\textwidth}
		\includegraphics[width=\textwidth]{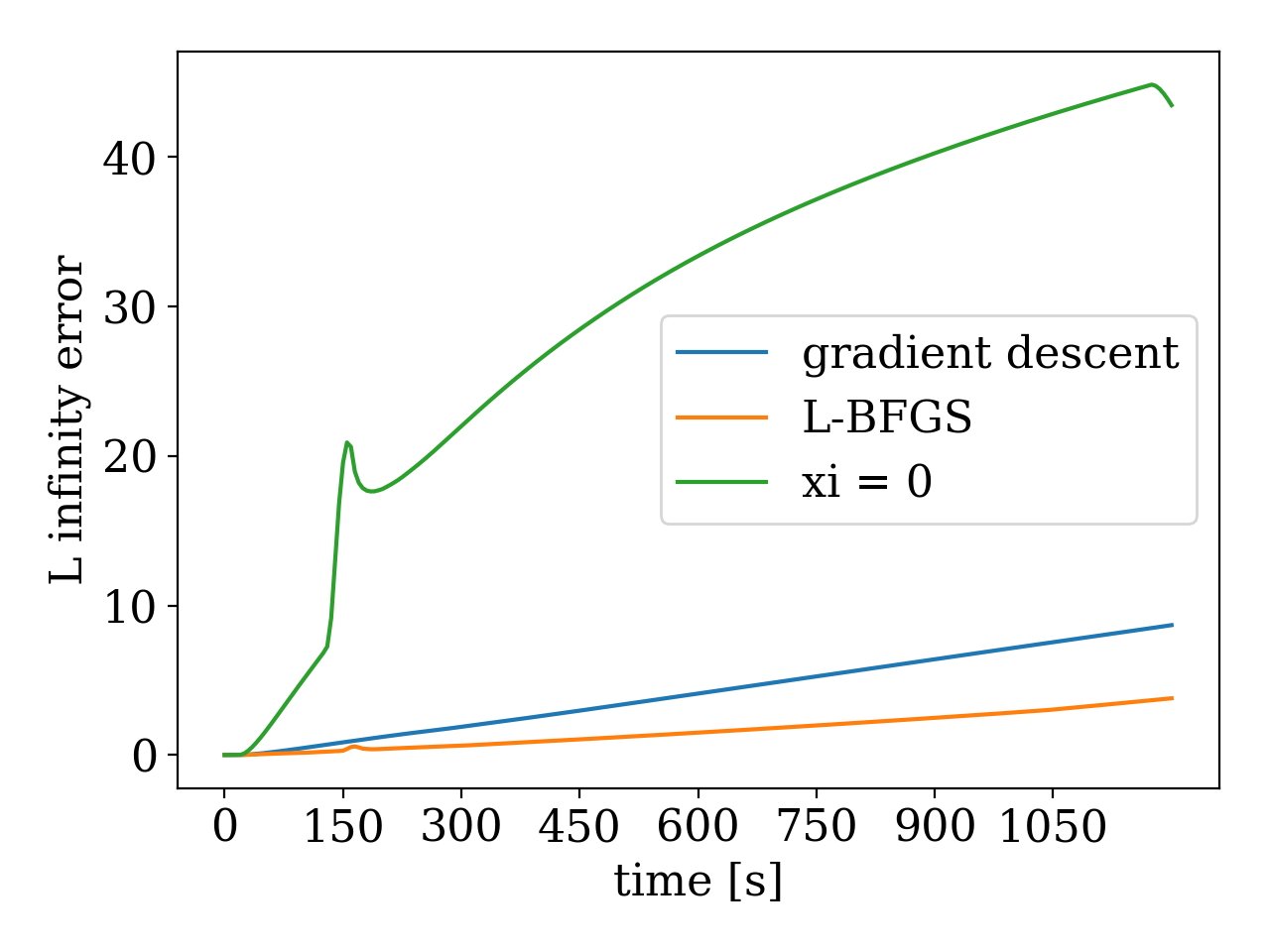}
		\caption{$L^\infty$-norm.}
	\end{subfigure}
	\hfil
	\begin{subfigure}{0.475\textwidth}
		\includegraphics[width=\textwidth]{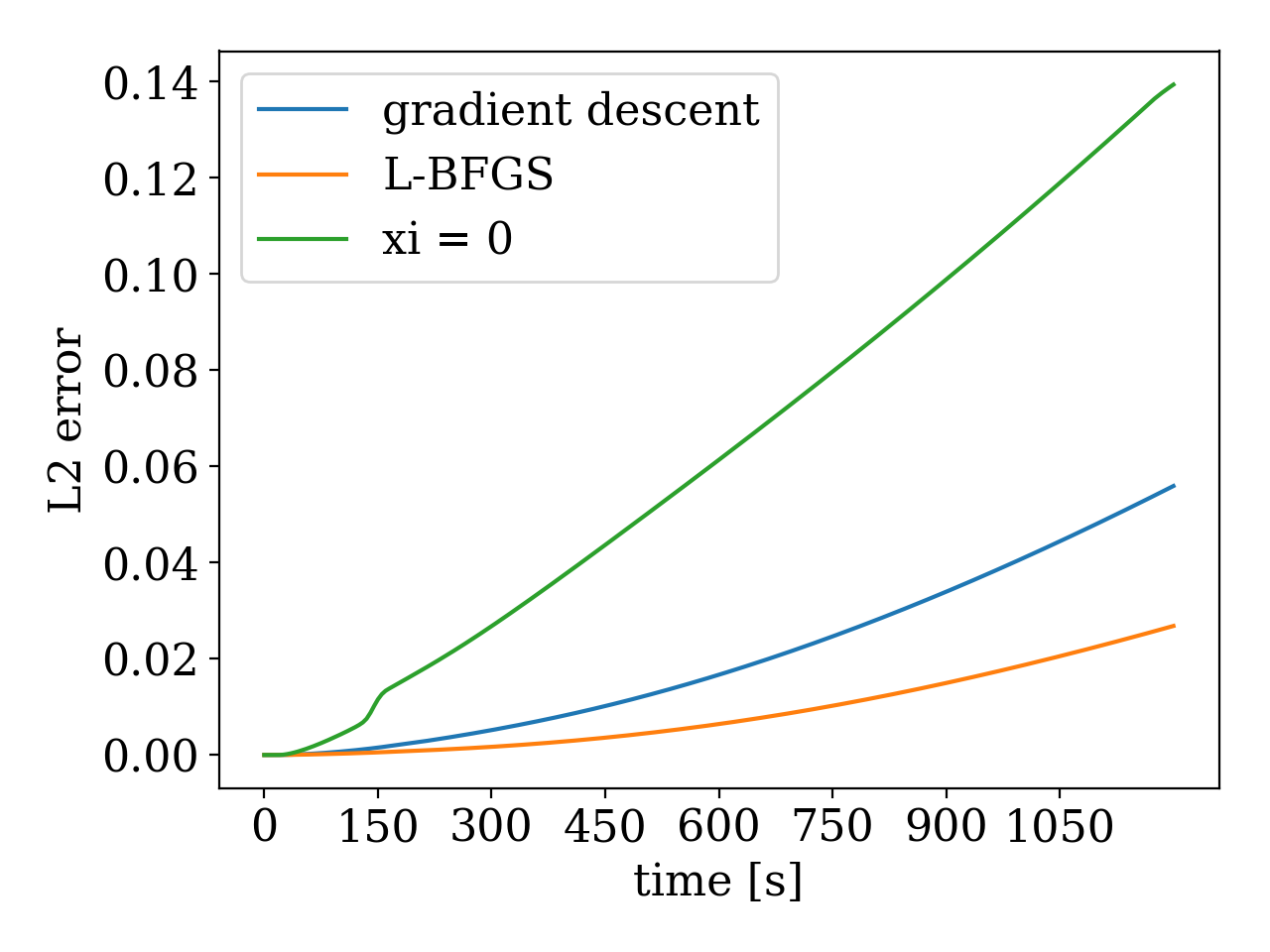}
		\caption{$L^2$-norm.}
	\end{subfigure}
	\caption{Evolution of the error in temperature over time for one measurement.}
	\label{fig:error_evolution_1}
\end{figure}

First, we note that the proposed parameter identification performs well for both algorithms. The identified perfusion rates in Figure~\ref{fig:identified_noiseless}  approximate the \qe{measured} perfusion rate well, at least close to the applicator. The positions of the blood vessels in the first row and next to the radiating boundary are the ones that are identified best, their position and shape closely resembles the measurement data. However, the identification becomes worse for the vessels located to the top and bottom of $\Gamma_\subrad$ and for the vessels further away from the applicator. In particular, the second row of vessels can only be seen very faintly for the gradient descent method, whereas the three middle vessels are still recognizable for the BFGS method, even though it underestimates the magnitude of the perfusion rate. This is due to the following reason: The influence of the perfusion rate is proportional to the temperature difference $\temperature - \temperature_\subblood$. However, the temperature is highest close to $\Gamma_\subrad$ and decays with increasing distance from there. Thus, the temperature difference at the \qe{outer} vessel locations is not significant and neither is the effect of the perfusion rate. Thus, our algorithm performs well by finding the significant blood vessels close to the applicator.

In Figure~\ref{fig:history_noiseless} we can further observe that the BFGS method outperforms the gradient descent algorithm as the values of the cost functions are always lower for the former and so are the values of the stationary measure. As the computational cost of the L-BFGS method is essentially the same as for the gradient descent algorithm, we can save valuable computational time by using the former, as it needs about half as many steps to reach a certain tolerance in the cost function compared to the latter.

The comparison of the simulated and measured physical quantities shown in Table~\ref{table:comparison_noiseless} emphasizes the fact, that the BFGS method exhibits better properties. Furthermore, comparing the errors of the simulated physical quantities to the ones we get for a vanishing perfusion rate we observe that the effect of blood perfusion is significant for the therapy planning as stated in \cite{perfusion_parameter}. We observe that the simulation results of our algorithms are significantly better than those of the simulation that does not consider the effect of the perfusion rate. Moreover, the errors generated by the BFGS method are only about half of those generated with the gradient descent method for $\temperature$ and $\radiation$, however, the error in tissue damage goes down dramatically by a factor of \num{3.5} in the $L^\infty$ norm and by a factor of \num{10} in the $L^2$ norm, underlining the superior behavior of the BFGS method. 

Finally, we also show the evolution of the error in temperature over time in Figure~\ref{fig:error_evolution_1} for both the $L^\infty(\Omega)$ and $L^2(\Omega)$ norm. We observe a similar picture as before, where the BFGS method outperforms the gradient descent method. Furthermore, we can see that even with only one measurement the methods produce comparatively low errors in the simulation, even though the results are \qe{extrapolated} into the future. In particular, we again observe that the effect of blood perfusion is significant here, as the error increases rapidly when considering $\xi = 0$, whereas the error stays comparatively small for the simulation with the identified perfusion rate.

\subsection*{Multiple Measurements}

\begin{figure}[!t]
	\centering
	\includegraphics[width=\textwidth]{legend}
	\begin{subfigure}[t]{0.475\textwidth}
		\includegraphics[width=\textwidth]{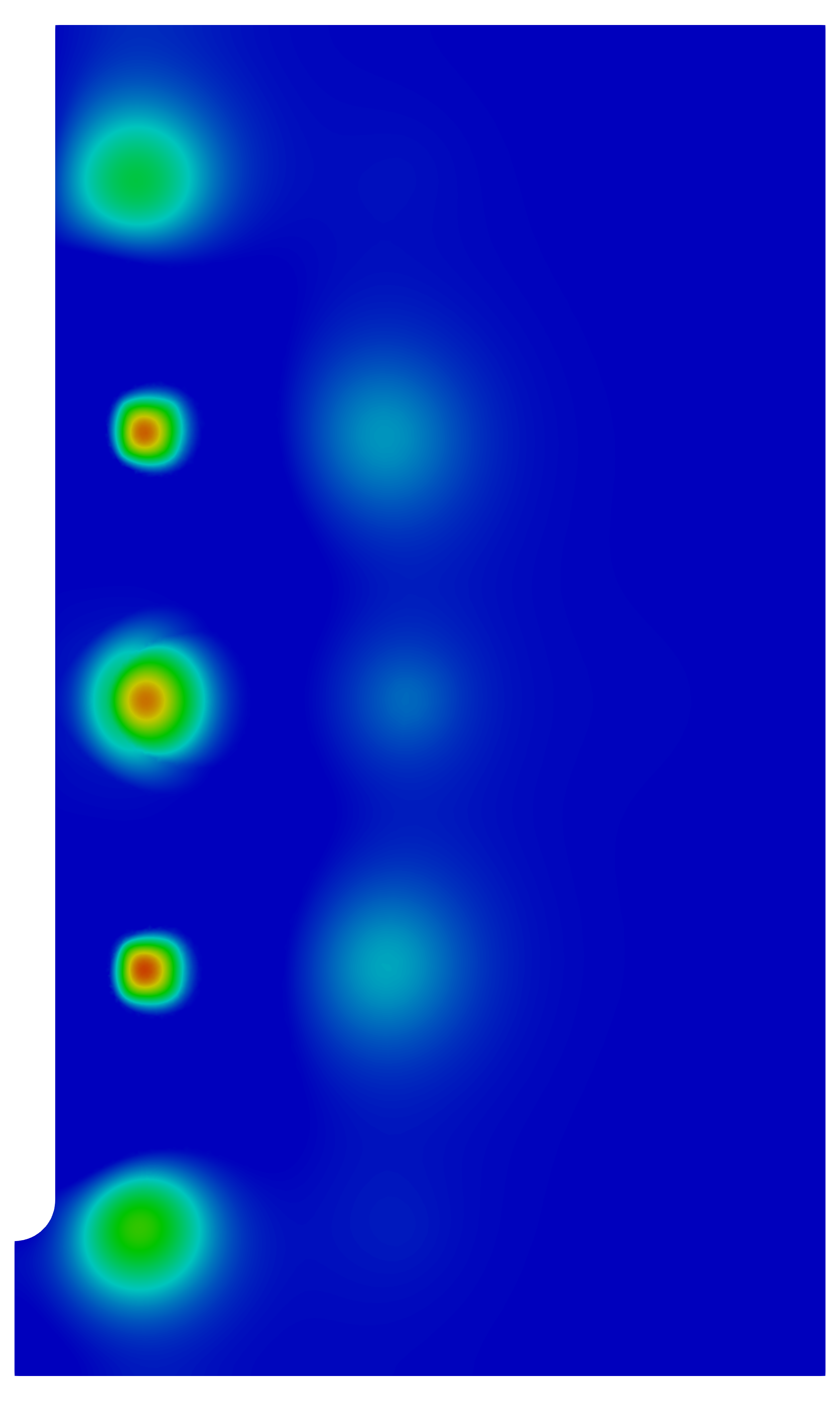}
		\caption{Projected gradient descent method.}
		\label{fig:identified_gradient_multiple}
	\end{subfigure}
	\hfil
	\begin{subfigure}[t]{0.475\textwidth}
		\includegraphics[width=\textwidth]{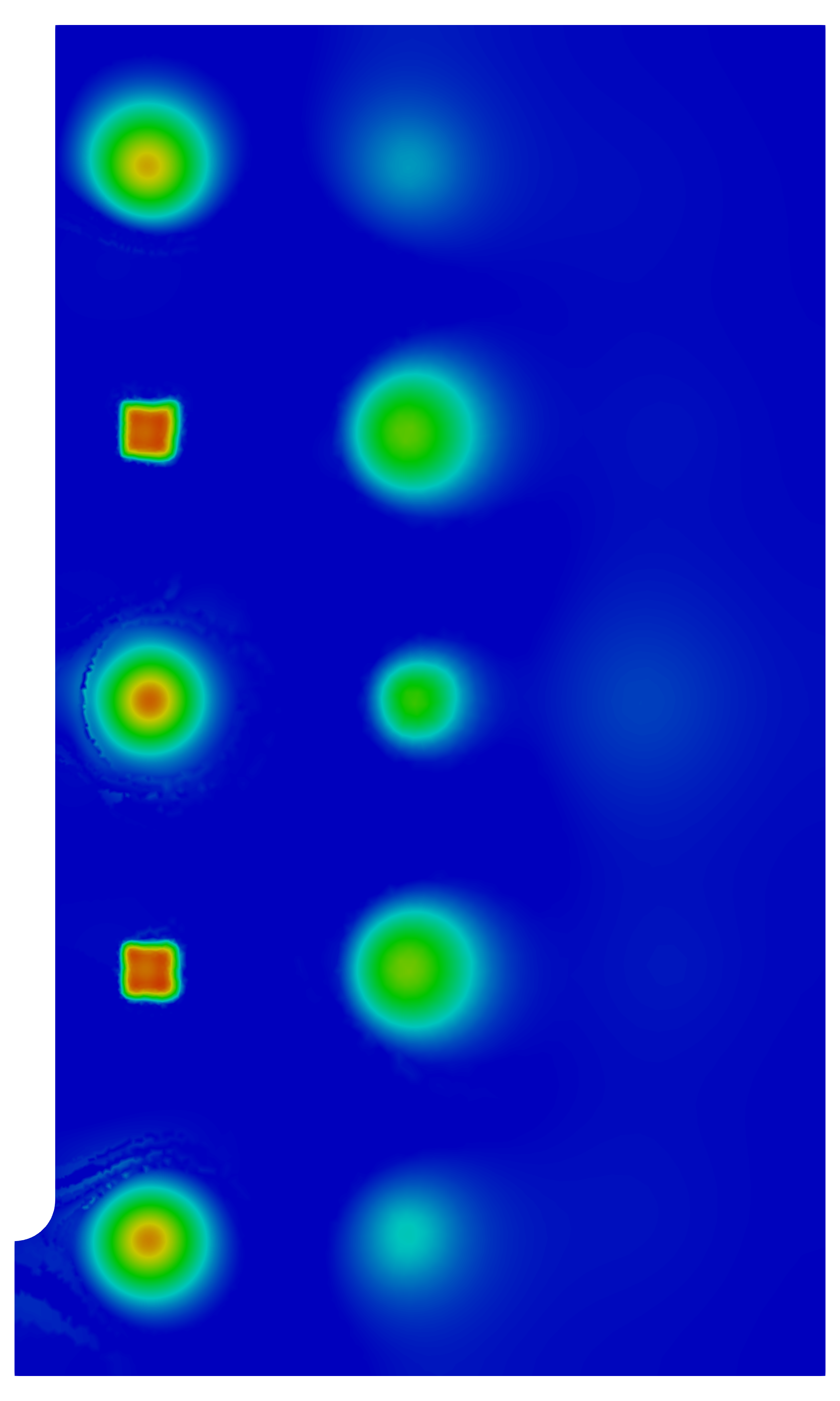}
		\caption{Projected L-BFGS method.}
		\label{fig:identified_bfgs_multiple}
	\end{subfigure}
	\caption{Identified perfusion rates for multiple measurements.}
	\label{fig:identified_noiseless_multiple}
\end{figure}

\begin{table}[!b]
	\centering
	\begin{tabular}{l r r }
		\toprule
		& $L^\infty$-error & $L^2$-error  \\
		\midrule
		gradient descent & & \\
		\midrule
		temperature $\temperature$ & \num{4.068} (\num{0.989} \%) & \num{0.466} (\num{0.372} \%) \\
		radiative energy $\radiation$ & \num{3819} (\num{1.053} \%) &  \num{5.914} (\num{0.097} \%) \\
		tissue damage $\dmgfun$ & \num{0.166} (\num{16.6} \%) & \num{0.002} (\num{1.934} \%) \\
		\midrule
		L-BFGS & & \\
		\midrule
		temperature $\temperature$ & \num{2.771} (\num{0.674} \%) & \num{0.253} (\num{0.2} \%) \\
		radiative energy $\radiation$ & \num{2443} (\num{0.674} \%) & \num{3.437} (\num{0.057} \%) \\
		tissue damage $\dmgfun$ & \num{0.105} (\num{10.46} \%) & \num{0.001} (\num{1.373} \%) \\
		\bottomrule
	\end{tabular}
	\caption{Comparison between the simulated and measured temperature for multiple measurements.}
	\label{table:comparison_noiseless_multiple}
\end{table}

Now, we also examine the quality of the identification for multiple measurements. We choose the same setting as before, but this time we assume to have temperature measurements at $\timehorizon_1 = \num{60}\ \si{\second}$, $\timehorizon_2 = \num{120}\ \si{\second}$ and $\timehorizon_3 = \num{180}\ \si{\second}$. The identified perfusion rates for this setting are shown in Figure~\ref{fig:identified_noiseless_multiple}, and the comparison to the measurement data is depicted in Table~\ref{table:comparison_noiseless_multiple}. 

Again, we observe that the L-BFGS method outperforms the gradient descent algorithm: In Figure~\ref{fig:identified_noiseless_multiple} we observe that it can identify more blood vessels and resolves them more accurately, e.g., the first column of vessels is clearly visible and well identified, the three vessels in the middle of the second column can be seen more clearly now, although their perfusion rate is still underestimated by the method. For the gradient descent algorithm only the middle vessels of the first column are identified accurately. The ones further away from $\Gamma_\subrad$ are better visible than before and the result resembles the one for the BFGS method with one measurement shown in Figure~\ref{fig:identified_bfgs}. 

The same is true for the comparison of the simulated and synthetic data (cf. Table~\ref{table:comparison_noiseless_multiple}). Here, the errors for the temperature with the three measurement gradient descent algorithm are very close to the errors obtained by the BFGS method with only one measurement. However, the errors for the radiative energy and tissue damage are lower. The BFGS method also improves its performance for three measurements compared to a single one, albeit only slightly for most errors. 

\begin{figure}[!t]
	\begin{subfigure}{0.475\textwidth}
		\includegraphics[width=\textwidth]{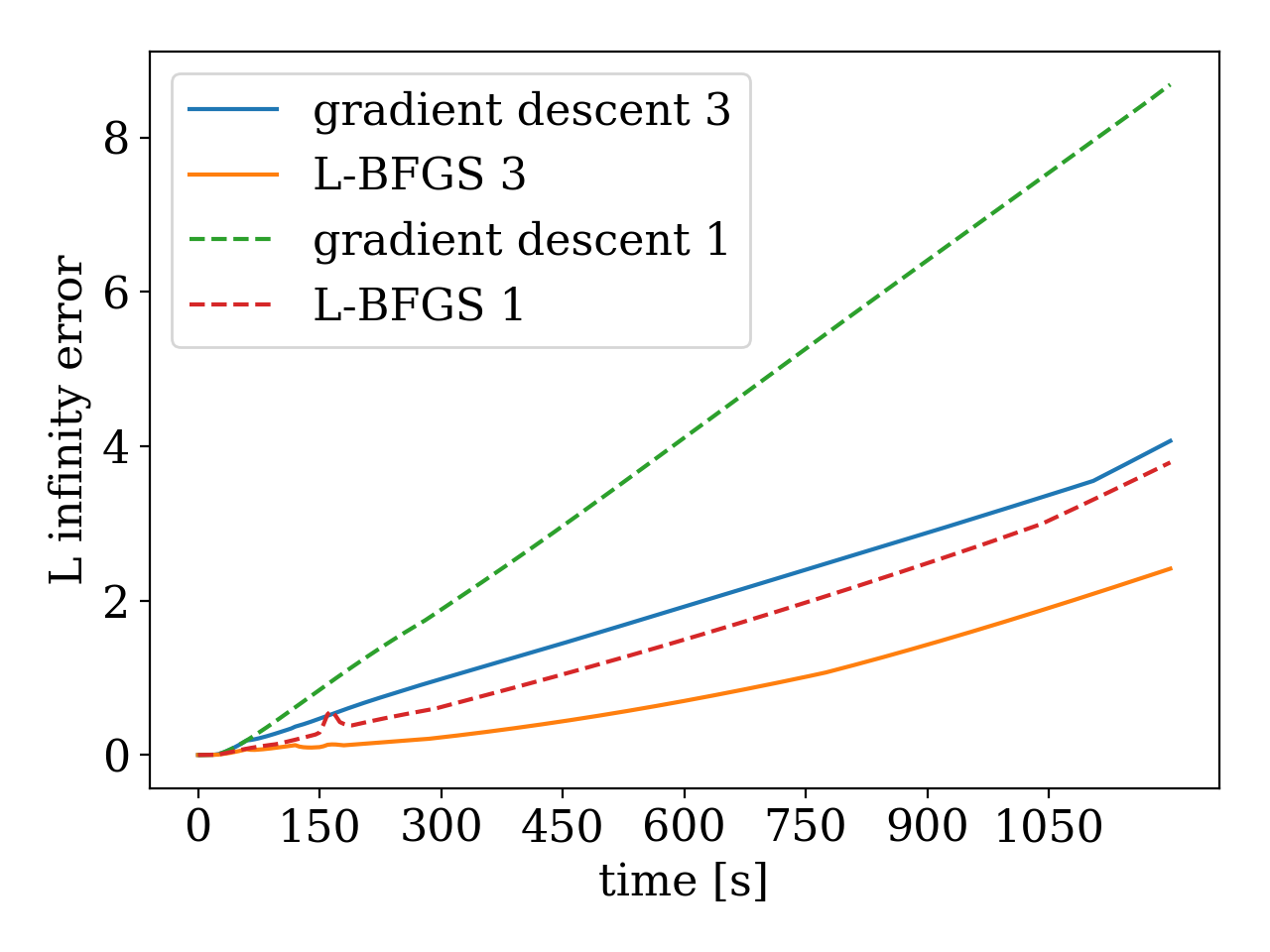}
		\caption{$L^\infty$-norm.}
	\end{subfigure}
	\hfil
	\begin{subfigure}{0.475\textwidth}
		\includegraphics[width=\textwidth]{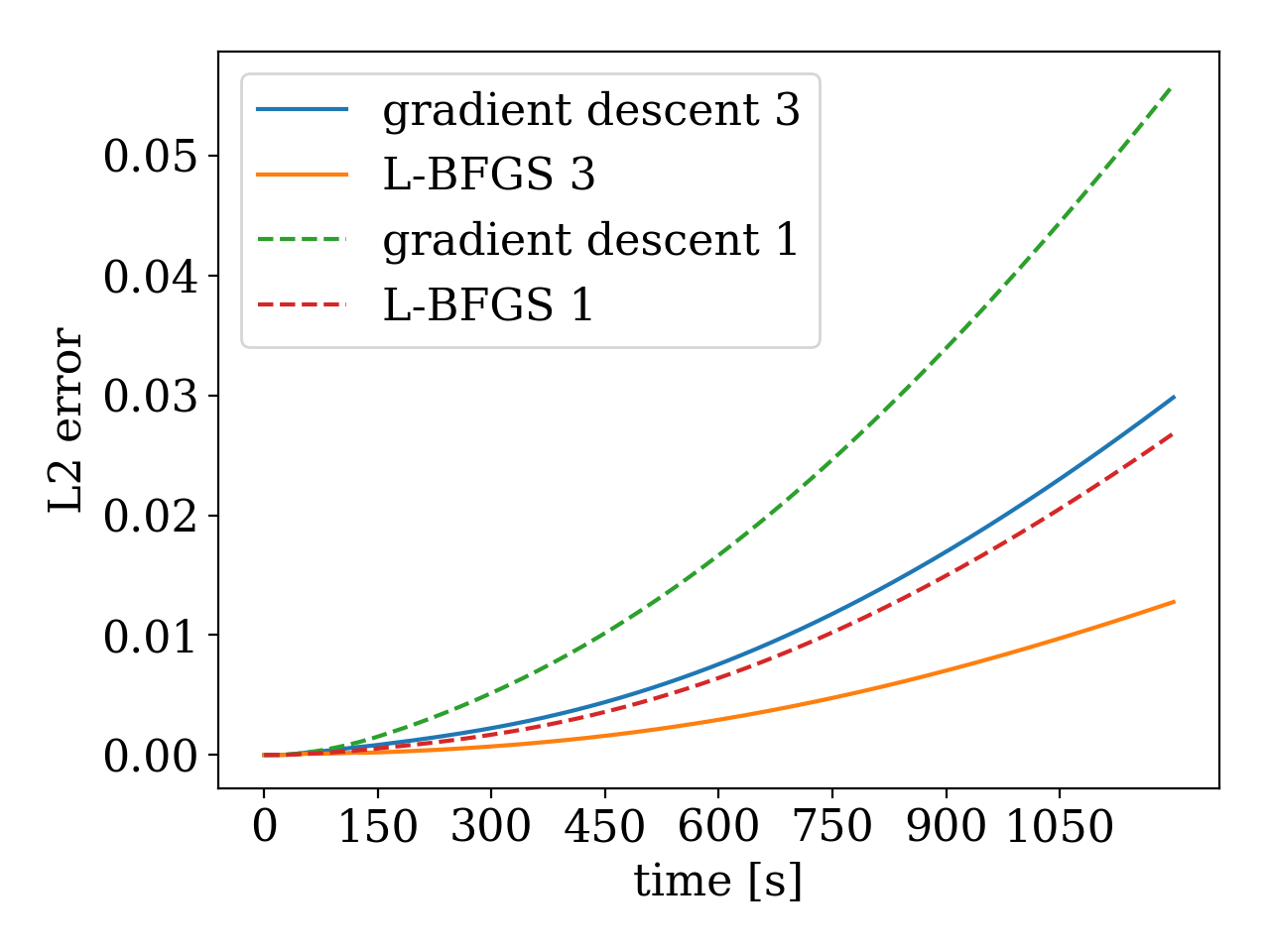}
		\caption{$L^2$-norm.}
	\end{subfigure}
	\caption{Evolution of the error in temperature over time for 3 measurements.}
	\label{fig:error_evolution_3}
\end{figure}

Again, we also show the evolution of the error in temperature over time for the whole simulation in Figure~\ref{fig:error_evolution_3}, this time in comparison to the results obtained with only one measurement. The results depicted here are similar to the ones of Figure~\ref{fig:error_evolution_1}. In particular, the BFGS method outperforms the gradient descent one. Moreover, we can see that the results obtained with the gradient descent method after three identifications are comparable to the ones of the BFGS method for only one identification process. It can also be seen that the error is smaller overall if we use three temperature measurements compared to using only one, and that it gives way better results even for late points in time during the therapy.

Summing up, we observe that a single measurement can be sufficient for our purposes if the identification is carried out with the BFGS method. Doing so saves a lot of computational time and is a stepping stone for the use of the method in an online therapy-planning and -monitoring tool.

\subsection{Noisy Model Problem}
\label{sec:numerics_noisy}

\begin{figure}[!b]
	\centering
	\includegraphics[width=\textwidth]{legend}
	\begin{subfigure}[t]{0.475\textwidth}
		\includegraphics[width=\textwidth]{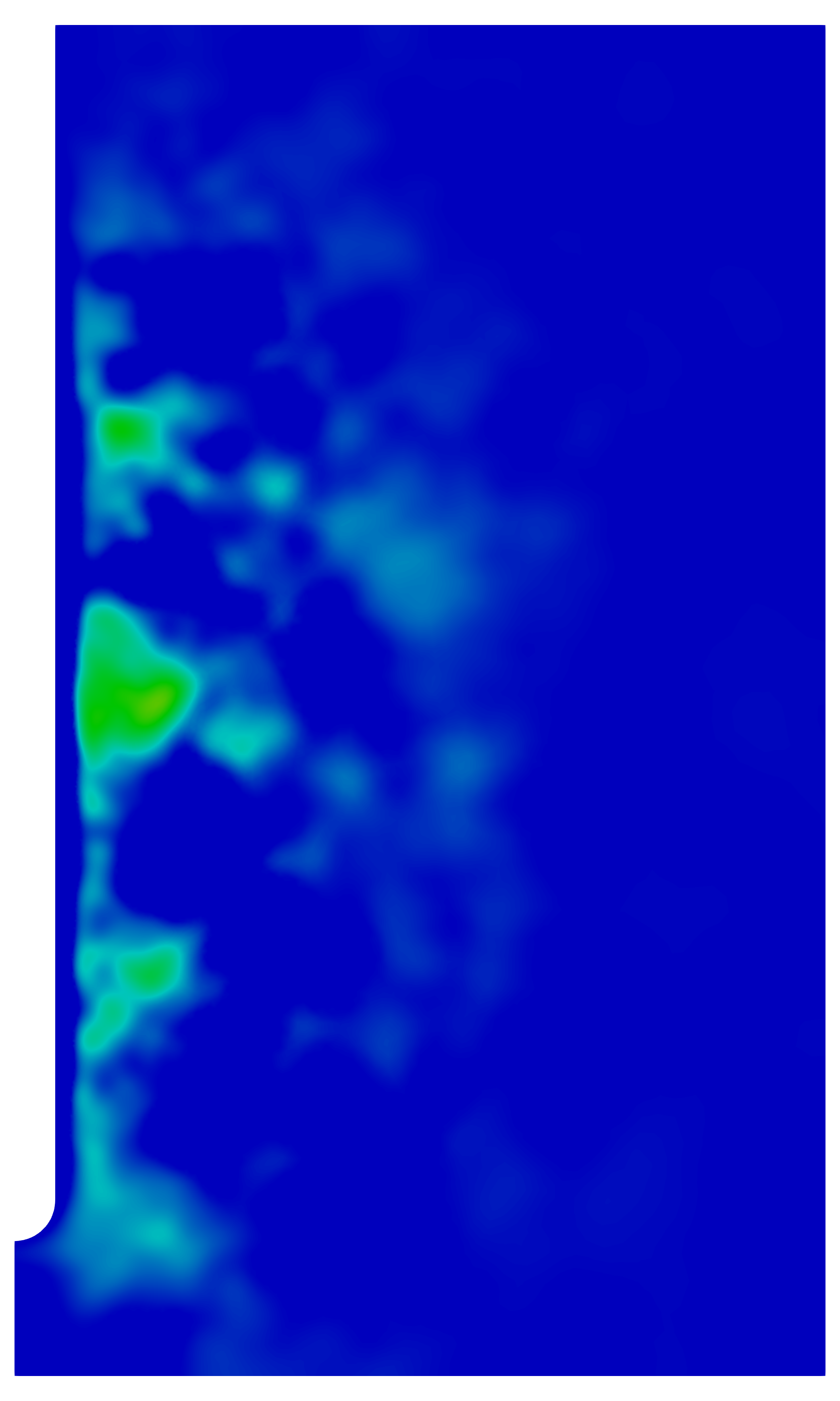}
		\caption{Projected gradient descent method.}
		\label{fig:identified_gradient_noisy}
	\end{subfigure}
	\hfil
	\begin{subfigure}[t]{0.475\textwidth}
		\includegraphics[width=\textwidth]{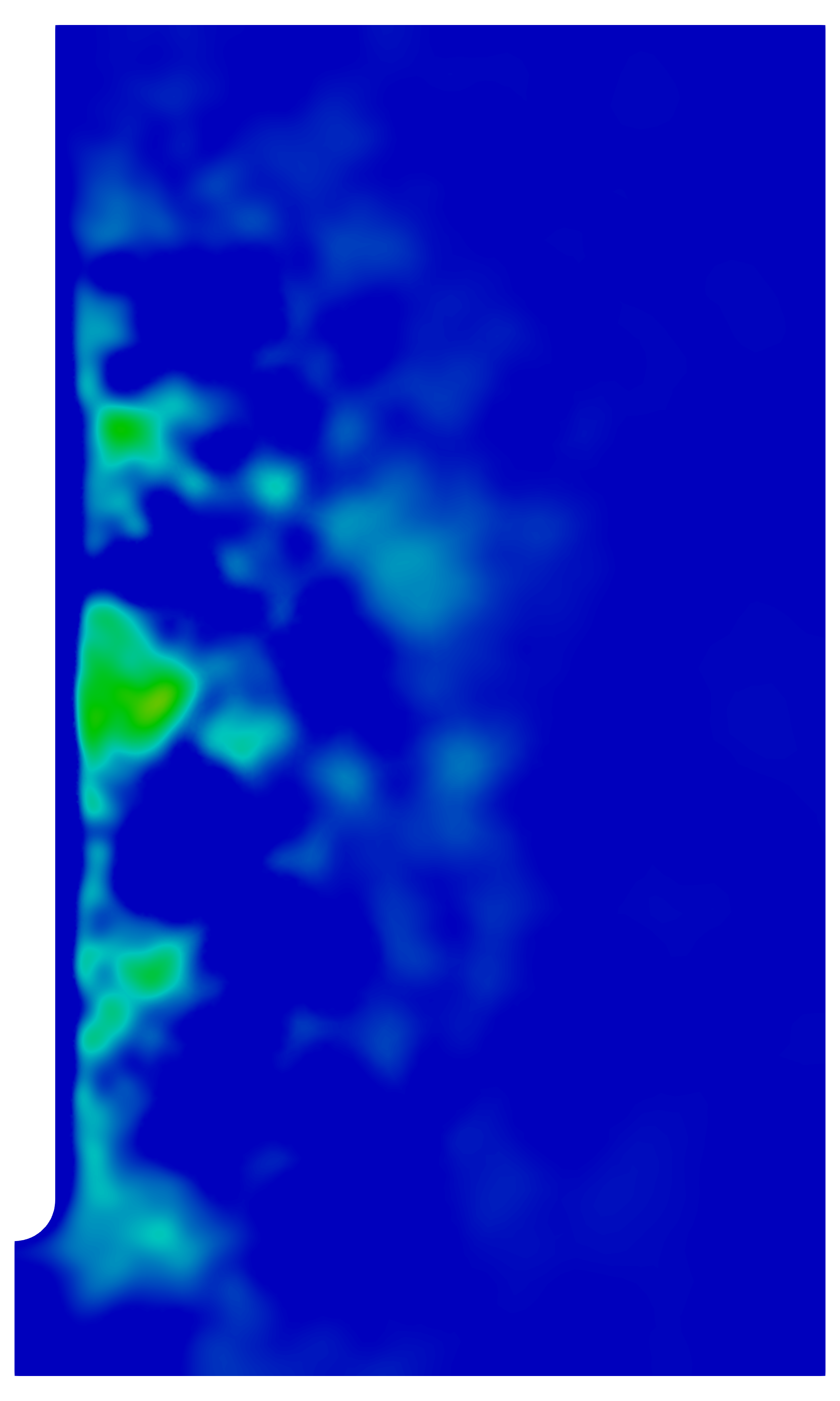}
		\caption{Projected L-BFGS method.}
		\label{fig:identified_bfgs_noisy}
	\end{subfigure}
	\caption{Identified perfusion rates for noisy measurements.}
	\label{fig:identified_noisy}
\end{figure}

Let us now consider the case of noisy measurement data: We choose the same desired perfusion rate as before (cf. Figure~\ref{fig:prescribed}) and generate the measurement by solving the state equations. Then, we add Gaussian noise with zero mean and a standard deviation of $\sigma = \num{2}~\si{\celsius}$. This is in accordance with \cite{thermometry_1} and \cite{thermometry_2} where the accuracy of temperature measurements using MR thermometry is reported to be \num{2.3}~\si{\celsius} and \num{1.3}~\si{\celsius}, respectively. To identify the perfusion rate for this problem, we first smooth the data with an isotropic linear diffusion process with end time \num{2e-7}, which is equivalent to convolving the noisy data with a Gaussian kernel that has a standard deviation of \num{6.32e-4}. Using a higher standard deviation leads to more smoothing but also introduces bias to the synthetic measurements such that the data is not compatible with the model anymore. For more details on linear diffusion processes for filtering noise we refer to, e.g., \cite{weickert}. As this problem is less regular we choose the regularization parameter as $\lambda = \num{2.5e-10}$. We compare the BFGS and gradient descent method in the same context as before: First, we have a single measurement of temperature at $\timehorizon = \num{60}\ \si{\second}$ and, second, we consider three measurements at $\timehorizon_1 = \num{60}\ \si{\second}$, $\timehorizon_2 = \num{120}\ \si{\second}$ and $\timehorizon_3 = \num{180}\ \si{\second}$. All other parameters are chosen as in Section~\ref{sec:numerics_noiseless}. 

\subsection*{A Single Measurement}

\begin{table}[!t]
	\centering
	\begin{tabular}{l r r }
		\toprule
		& $L^\infty$-error & $L^2$-error  \\
		\midrule
		gradient descent & & \\
		\midrule
		temperature $\temperature$ & \num{8.805} (\num{2.141} \%) & \num{0.737} (\num{0.588} \%) \\
		radiative energy $\radiation$ & \num{111237} (\num{30.67} \%) &  \num{159} (\num{2.62} \%) \\
		tissue damage $\dmgfun$ & \num{0.973} (\num{97.27} \%) & \num{0.016} (\num{17.879} \%) \\
		\midrule
		L-BFGS & & \\
		\midrule
		temperature $\temperature$ & \num{8.876} (\num{2.158} \%) & \num{0.712} (\num{0.568} \%) \\
		radiative energy $\radiation$ & \num{112285} (\num{31} \%) & \num{160} (\num{2.64} \%) \\
		tissue damage $\dmgfun$ & \num{0.974} (\num{97.4} \%) & \num{0.017} (\num{18.3} \%) \\
		\bottomrule
	\end{tabular}
	\caption{Comparison between the simulated and measured temperature for a single measurement.}
	\label{table:comparison_noisy}
\end{table}

\begin{figure}[!b]
	\begin{subfigure}{0.475\textwidth}
		\includegraphics[width=\textwidth]{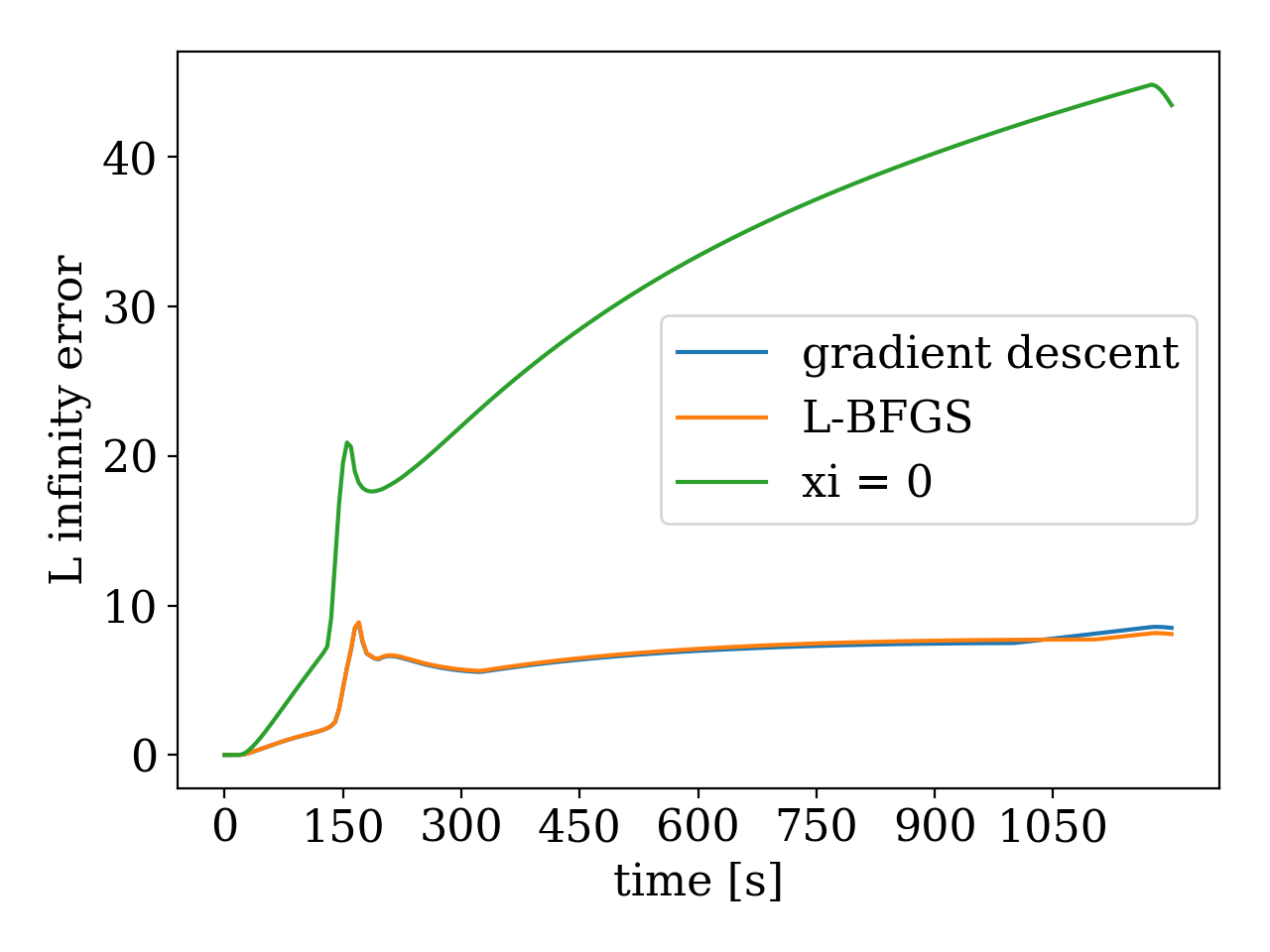}
		\caption{$L^\infty$-norm.}
	\end{subfigure}
	\hfil
	\begin{subfigure}{0.475\textwidth}
		\includegraphics[width=\textwidth]{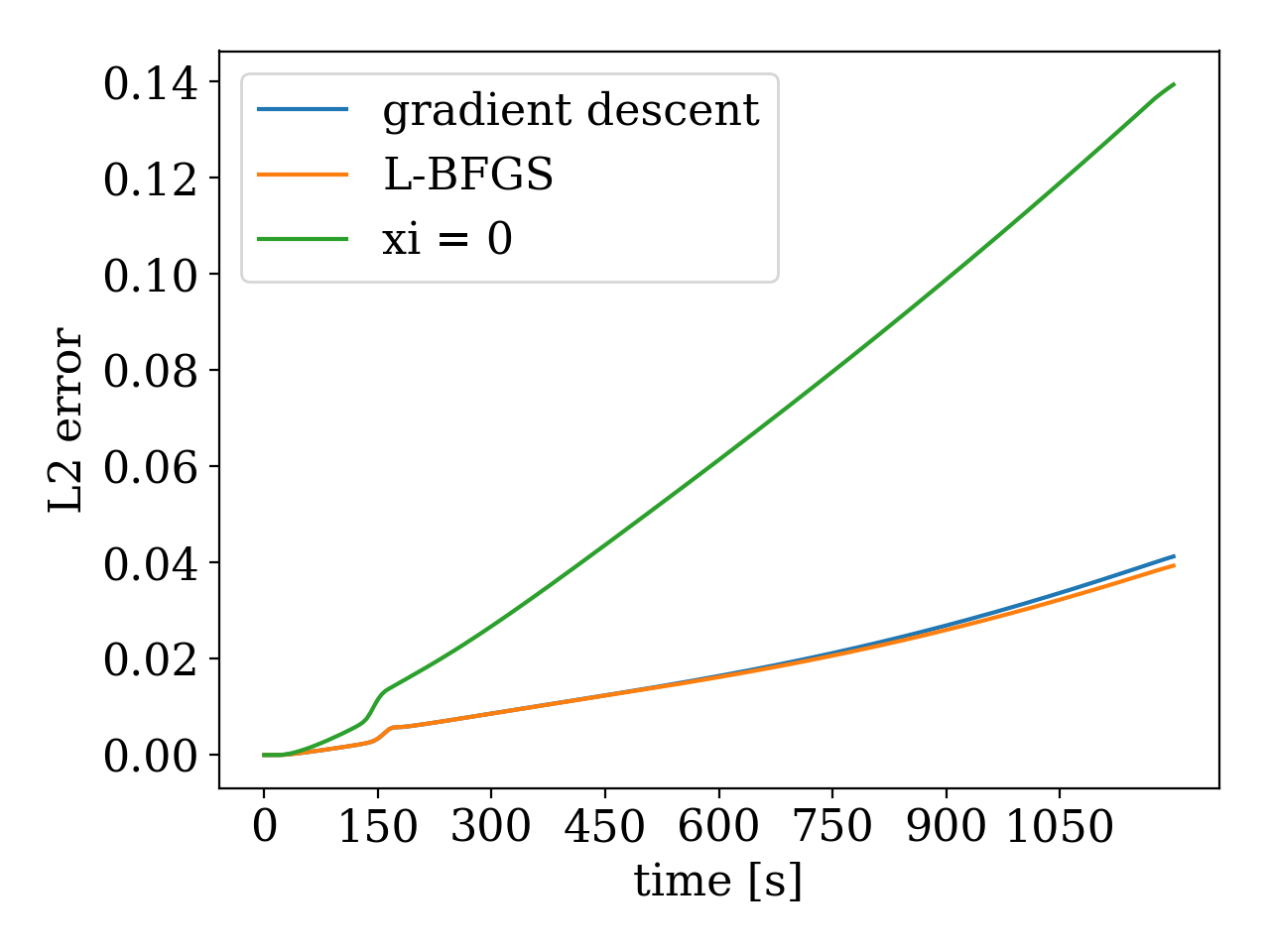}
		\caption{$L^2$-norm.}
	\end{subfigure}
	\caption{Evolution of the error in temperature over time for one noisy measurement.}
	\label{fig:error_evolution_1_noisy}
\end{figure}

In the case of a single measurement at $\timehorizon = \num{60}\ \si{\second}$ we see that there is not much of a difference between the gradient descent method and the L-BFGS algorithm. The computed perfusion rates, depicted in Figure~\ref{fig:identified_noisy}, look nearly identical, and so do the errors, as can be seen in Table~\ref{table:comparison_noisy} and Figure~\ref{fig:error_evolution_1_noisy}. As is to be expected, the results obtained from the noisy measurements are worse than the ones we got in Section~\ref{sec:numerics_noiseless}, where we did only consider noiseless data. In particular, the computed perfusion rates are worse compared to the previous cases and one can hardly distinguish the isolated vessels anymore. Instead, the vessels are more \qe{smeared out}. However, we can still observe that we have three peaks in the perfusion rate, corresponding to the vessels in the middle of the first column for the desired perfusion rate. Additionally, we can observe that the middle one of these three shows a higher perfusion rate, as it is the case for the synthetic blood perfusion (cf. Figure~\ref{fig:prescribed}). We also observe that both methods underestimate the perfusion rate by about \num{20} \%. This is due to the effect of the regularization term that penalizes large values of $\perfusion$. 

The quality of approximation for the whole therapy, as shown in Table~\ref{table:comparison_noisy}, is comparable for both methods, neither one of them performs significantly better than the other. The same is true when investigating the evolution of the error in temperature over time (cf. Figure~\ref{fig:error_evolution_1_noisy}), where there is barely any difference visible between both algorithms.

\subsection*{Multiple Measurements}

\begin{figure}[!t]
	\centering
	\includegraphics[width=\textwidth]{legend}
	\begin{subfigure}[t]{0.475\textwidth}
		\includegraphics[width=\textwidth]{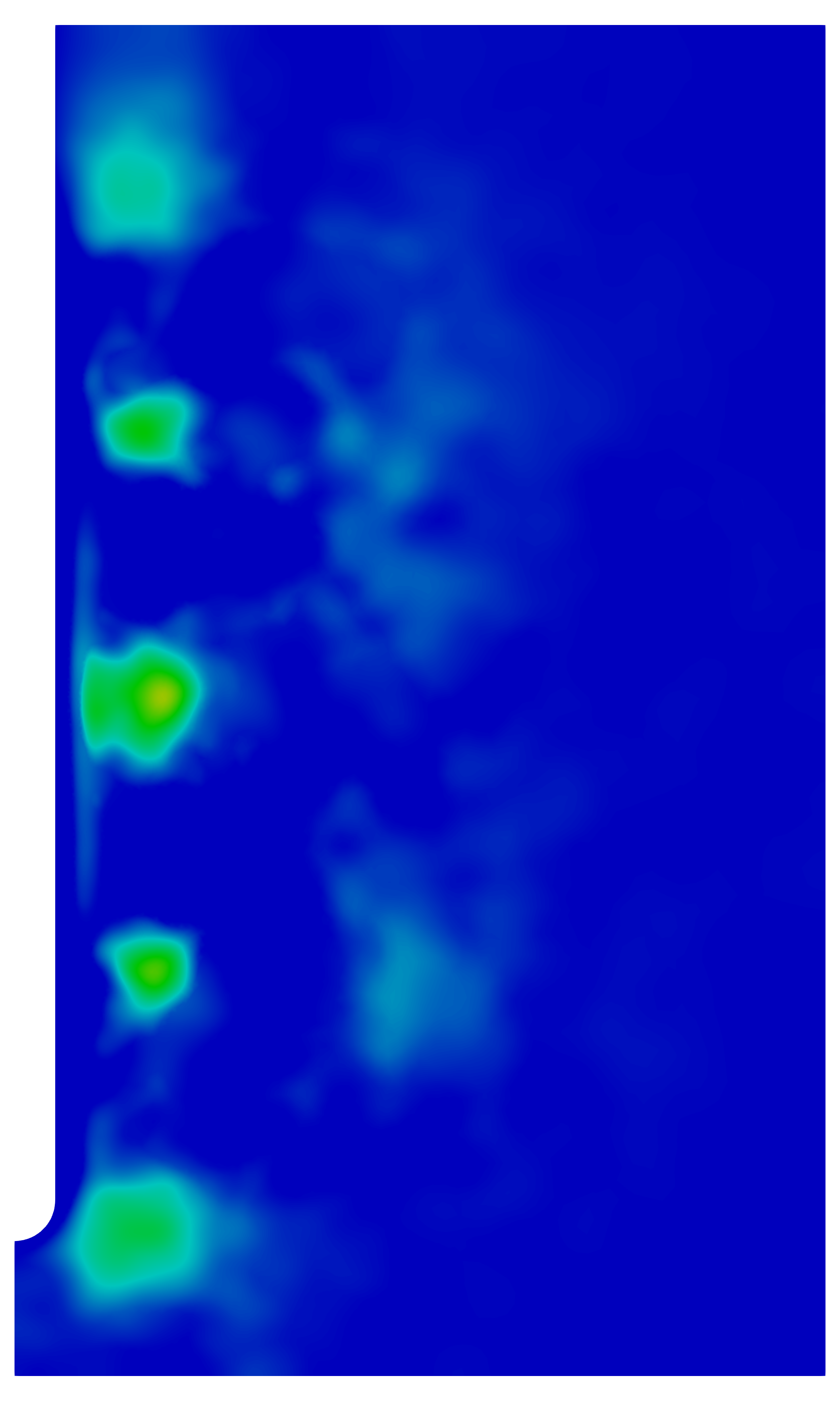}
		\caption{Projected gradient descent method.}
		\label{fig:identified_gradient_multiple_noisy}
	\end{subfigure}
	\hfil
	\begin{subfigure}[t]{0.475\textwidth}
		\includegraphics[width=\textwidth]{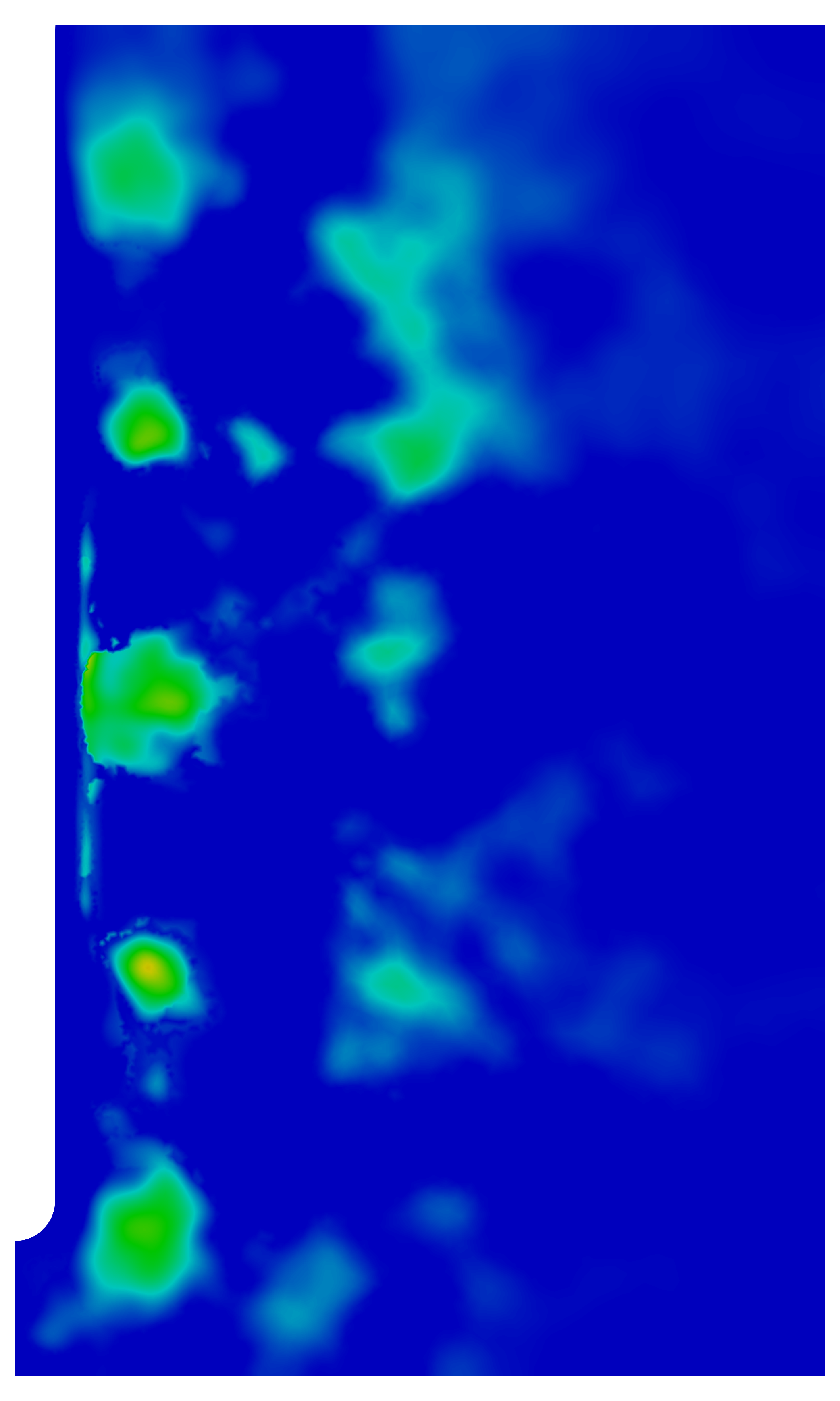}
		\caption{Projected L-BFGS method.}
		\label{fig:identified_bfgs_multiple_noisy}
	\end{subfigure}
	\caption{Identified perfusion rates for multiple noisy measurements.}
	\label{fig:identified_noisy_multiple}
\end{figure}

\begin{table}[!b]
	\centering
	\begin{tabular}{l r r }
		\toprule
		& $L^\infty$-error & $L^2$-error  \\
		\midrule
		gradient descent & & \\
		\midrule
		temperature $\temperature$ & \num{5.148} (\num{1.25} \%) & \num{0.51} (\num{0.41} \%) \\
		radiative energy $\radiation$ & \num{38413} (\num{10.6} \%) &  \num{50} (\num{0.823} \%) \\
		tissue damage $\dmgfun$ & \num{0.589} (\num{58.9} \%) & \num{0.006} (\num{7.16} \%) \\
		\midrule
		L-BFGS & & \\
		\midrule
		temperature $\temperature$ & \num{4.374} (\num{1.06} \%) & \num{0.4} (\num{0.32} \%) \\
		radiative energy $\radiation$ & \num{40890} (\num{11.27} \%) & \num{50.15} (\num{0.825} \%) \\
		tissue damage $\dmgfun$ & \num{0.62} (\num{62} \%) & \num{0.006} (\num{6.93} \%) \\
		\bottomrule
	\end{tabular}
	\caption{Comparison between the simulated and measured temperature for multiple noisy measurements.}
	\label{table:comparison_noisy_multiple}
\end{table}

For the case of three noisy measurements, both the gradient descent algorithm and the L-BFGS method perform significantly better and there are more differences between the methods. In particular, the obtained perfusion rates are much closer to the synthetic one as Figure~\ref{fig:identified_noisy_multiple} shows. Here, we observe that both methods can resolve the individual peaks of the blood vessels placed in the first row, now even all five of them are well visible. Again, we observe that the algorithms underestimate the perfusion rate, however, the error is smaller for the L-BFGS method. We can also see that there are some artifacts contained in the computed perfusion rates. We observe that noise is added to the regions behind the first column, especially in the L-BFGS method.

The quality of approximation for the whole therapy also improves, as can be seen in Table~\ref{table:comparison_noisy_multiple}. Again, both methods show significant improvements, in particular for the error in temperature, which is about halved in the $L^\infty$-norm, and in tissue damage, which decreases nearly by two-thirds in the $L^2$-norm. The evolution of the temperature error also confirms these findings, as the results are substantially better as for only one measurement. Additionally, we can see that the L-BFGS method now also performs better than the gradient descent algorithm.

\begin{figure}[!t]
	\begin{subfigure}{0.475\textwidth}
		\includegraphics[width=\textwidth]{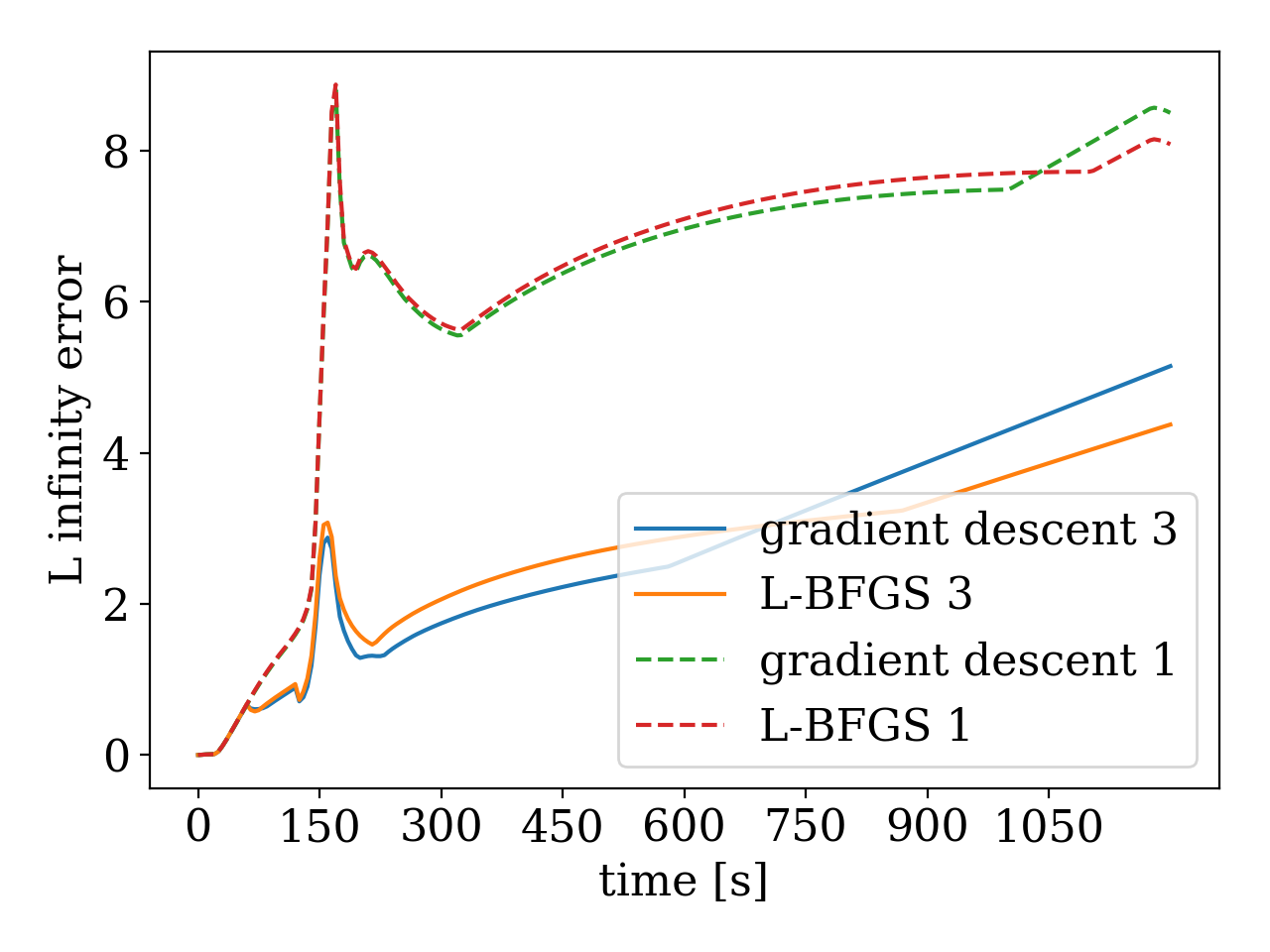}
		\caption{$L^\infty$-norm.}
	\end{subfigure}
	\hfil
	\begin{subfigure}{0.475\textwidth}
		\includegraphics[width=\textwidth]{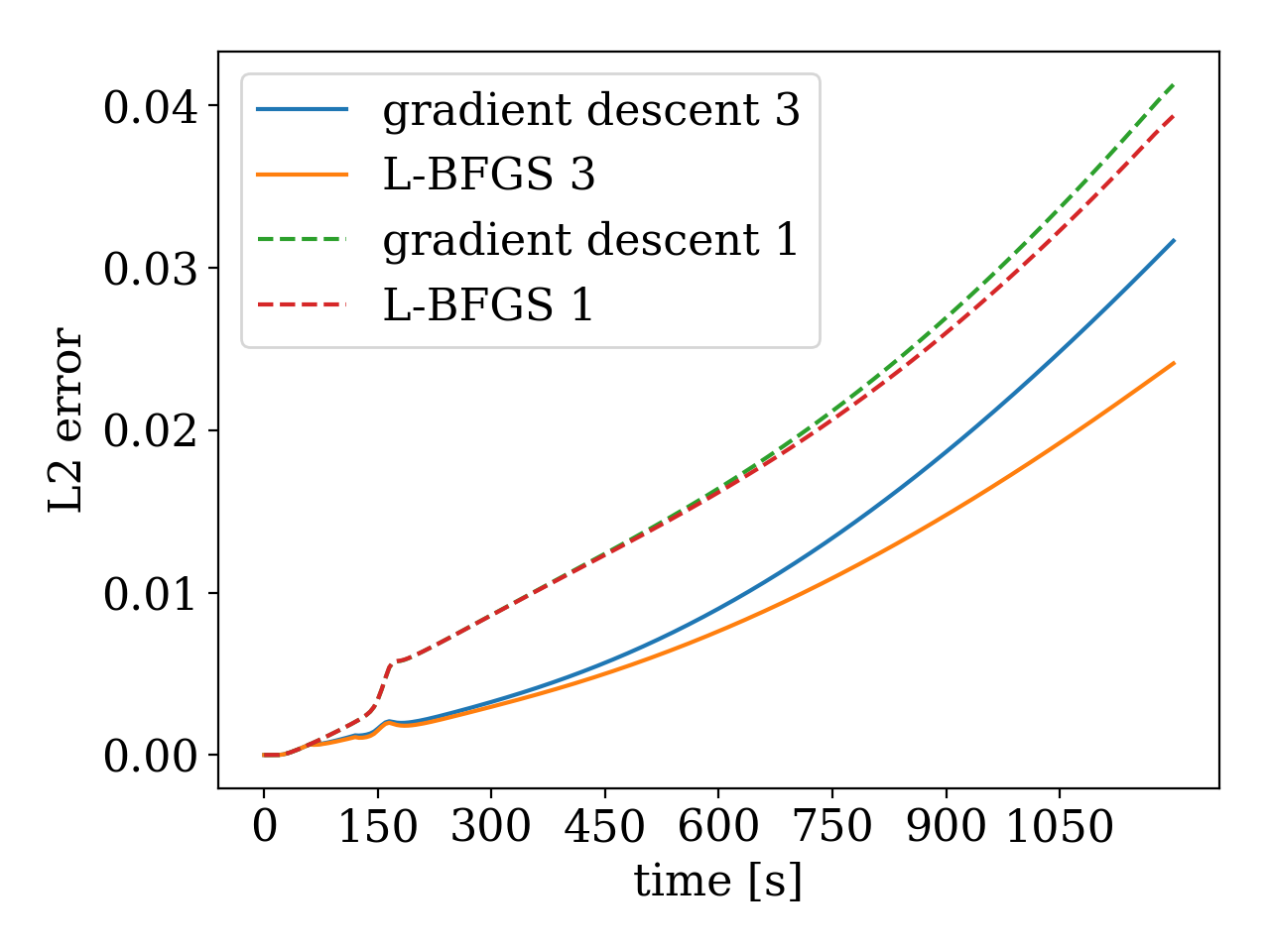}
		\caption{$L^2$-norm.}
	\end{subfigure}
	\caption{Evolution of the error in temperature over time for 3 noisy measurements.}
	\label{fig:error_evolution_3_noisy}
\end{figure}

\section{Conclusions}
\label{sec:conclusion}
We have demonstrated that the proposed parameter identification approach based on techniques from PDE-constrained optimization can identify the blood perfusion rate in the relevant region around the applicator. This was done using synthetic measurements with and without artificial noise. Making use of three instead of one subsequent measurements has notably improved the accuracy. The L-BFGS method uses significantly less iterations to converge to an acceptable solution than the gradient descent method, while the time per iteration is comparable for both methods. The next step will be to test the parameter identification with real MR thermometry data obtained from ex-vivo experiments with artificial blood vessels.


\begin{backmatter}

\section*{Competing interests}
  The authors declare that they have no competing interests.
  

\section*{Author's contributions}
    All authors read and approved the final manuscript.

\section*{Acknowledgements}

The authors acknowledge the financial support by the Federal Ministry of Education and Research of Germany in the framework of the project \textit{proMT: Prognostische modellbasierte online MR-Thermometrie bei minimalinvasiver Thermoablation zur Behandlung von Lebertumoren} (grant no's: 05M16AMA and 05M16UKE). 


\bibliographystyle{bmc-mathphys} 
\bibliography{references.bib}      

\end{backmatter}
\end{document}